
\magnification=1200

\input amstex
\documentstyle{amsppt}
\define\CDdashright#1#2{&\,\mathop{\dashrightarrow}\limits^{#1}_{#2}\,&}
\define\CDdashleft#1#2{&\,\mathop{\dashleftarrow}\limits^{#1}_{#2}\,&}

\def\P{\Bbb P}

\def\Til#1{{\widetilde{#1}}}
\def\codim{\text{codim}}


\input xypic

\def\cmp{{c_{\text{MP}}}}

\def\ddi#1#2{\frac{\partial {#1}}{\partial x_{#2}}}
\def\ddti#1#2{\frac{\partial \Til{#1}}{\partial \tilde x_{#2}}}
\def\ddfi#1{\frac{\partial F}{\partial x_{#1}}}
\def\ddft#1{\frac{\partial F}{\partial \tilde x_{#1}}}
\def\red#1{{#1}_{\text{red}}}

\def\NB{{\text{\sl NB}}}

\def\siltable#1.{
\vbox{\tabskip=0pt \offinterlineskip
\halign to360pt{\strut##& ##\tabskip=1em plus2em&
  \hfil##\hfil& \vrule##&
  \hfil##\hfil& \vrule##&
  \hfil##\hfil& \vrule##&
  \hfil##\hfil& \vrule\thinspace\vrule##&
  \hfil##\hfil& \vrule\thinspace\vrule##&
  \hfil##\hfil& ##\tabskip=0pt\cr
#1}}}
\def\bdll{\Cal P^1_M\Cal L}
\def\bdlr{\Cal P^1_{\Til M}\Cal L}
\def\Grass{\text{Grass}}

\Monograph

\topmatter
\title Chern classes for singular hypersurfaces\endtitle
\author Paolo Aluffi${}^1$\endauthor
\date June 1996\enddate
\affil Florida State University\endaffil
\address Tallahassee, FL 32306\endaddress
\email aluffi\@ math.fsu.edu\endemail
\abstract We prove a simple formula for MacPherson's Chern class of
{\it hypersurfaces\/} in nonsingular varieties. The result highlights
the relation between MacPherson's class and other definitions of
homology Chern classes of singular varieties, such as Mather's Chern
class and the class introduced by W. Fulton in \cite{Fulton},
4.2.6.\endabstract
\endtopmatter
\document
{\eightpoint\toc
\widestnumber\head{\S4.}
\widestnumber\subhead{\S3.3.}
\head \S0. Introduction\page{2}\endhead
\head \S1. Statements of the result\page{4}\endhead
\subhead \S1.1. MacPherson's Chern class and Segre classes\page{4}\endsubhead
\subhead \S1.2. MacPherson's Chern class and Fulton's Chern
class\page{5}\endsubhead
\subhead \S1.3. MacPherson's Chern class and Mather's Chern
class\page{6}\endsubhead
\subhead \S1.4. Notational device\page{8}\endsubhead
\subhead \S1.5. MacPherson's Chern class and $\mu$@-classes\page{9}\endsubhead
\head \S2. The proof: preliminaries and divisors with normal
crossings\page{10}\endhead
\subhead \S2.1. $c_*(X)=c_*(\red X)$\page{11}\endsubhead
\subhead \S2.2. Divisors with normal crossings: proof\page{13}\endsubhead
\head \S3. The proof: behavior under blow@-ups\page{16}\endhead
\subhead \S3.1. (3) in terms of classes in $\P\bdll$, $\P\bdlr$\page{16}\endsubhead
\subhead \S3.2. (3) in terms of classes in $\P(\bdll\oplus\bdlr)$ over
$\Til M$\page{18}\endsubhead
\subhead \S3.3. (3) in terms of classes in $\P(\bdll\oplus\bdlr)$ over
$BL$\page{20}\endsubhead
\subhead \S3.4. The graph construction\page{21}\endsubhead
\subhead \S3.5. Coordinate set@-up in $BL$\page{22}\endsubhead
\subhead \S3.6. Three lemmas\page{26}\endsubhead
\subhead \S3.7. Computing $Z_\infty$\page{28}\endsubhead
\subhead \S3.8. End of the proof of (3)\page{31}\endsubhead
\head \S4. Remarks and applications\page{33}\endhead 
\subhead \S4.1. $\mu$@-class and Parusi\'nski's Milnor number\page{33}\endsubhead
\subhead \S4.2. Blowing up $\mu$@-classes\page{34}\endsubhead
\subhead \S4.3. Contact of two hypersurfaces\page{35}\endsubhead
\subhead \S4.4. A geometric application\page{37}\endsubhead
\endtoc}

\document
\footnote[]{${}^1$Supported in part by NSF grant DMS-9500843}

\newpage


\head \S0. Introduction\endhead

There are several candidates for a notion of homology `Chern classes'
of a (possibly singular) algebraic variety $X$, all agreeing with the
classes obtained as duals of the Chern classes of the tangent bundle of
$X$ if $X$ is nonsingular. Robert MacPherson (\cite{MacPherson})
introduced one such notion (which we will denote $\cmp(X)$) in proving a
conjecture (attributed to P.~Deligne and A. Grothendieck) prescribing a
functorial set@-up which makes it extremely well@-behaved. The purpose of
this note is to prove a simple formula for MacPherson's class of a {\it
hypersurface} of a nonsingular variety.

A summary of the context and essential definitions is given in \S1,
where we will state the result in several different forms. In this
introduction we will state a version of the result which puts in the same
framework MacPherson's class of a hypersurface $X$ together with two
other notions of homology Chern classes of $X$. Suppose $X$ is a
hypersurface of a nonsingular variety $M$ over an algebraically closed
field of characteristic 0. By the `singular subscheme' of $X$ we mean
the subscheme $Y$ of $X$ locally defined by the partial derivatives of
an equation for $X$. Now consider the blow@-up of $M$ along $Y$; this is
a variety mapping to $M$, and with several interesting divisors:

---the exceptional divisor, which we will denote $\Cal Y$;

---the scheme@-theoretic inverse image of $X$, still denoted $X$;

---the proper transform of $X$, denoted $\NB$.

It is well@-known that $\NB$ is in fact a realization of the Nash
blow@-up of $X$. If $\pi$ denotes the map to $M$ (and its
restrictions), we may use this set@-up to construct three classes living
in the Chow group $A_*X$ of $X$:
$$\align
c_F(X)=c(TM)\cap &\pi_*\frac{[X]}{1+X} \tag{*}\\
c_M(X)=c(TM)\cap &\pi_*\frac{[\NB]}{1+X-\Cal Y}\tag{**}\\
c_*(X)=c(TM)\cap &\pi_*\frac{[X]-[\Cal Y]}{1+X-\Cal Y}\tag{***}
\endalign$$
where $c(TM)$ denotes the total Chern class of the tangent bundle of $M$.

If $X$ is nonsingular, then $Y=\emptyset$ and these
classes coincide vacuously; they are seen at once to give then the
homology total Chern class of $X$
$$c(TX)\cap [X]$$
The class $c_F(X)$ listed above is {\it Fulton's Chern class\/} of $X$ (see
\S1.2), defined in \cite{Fulton}, Example 4.2.6: indeed, for hypersurfaces
the latter is simply the class of the virtual tangent bundle of $X$, which
(*) computes. The class (**) is easily seen to
compute {\it Mather's Chern class\/} $c_M(X)$ of $X$ (see \S1.3). The
main result of this note is:
\proclaim{Theorem} $c_*(X)$ equals MacPherson's Chern class of
$X$.\endproclaim

It seems very remarkable to us that MacPherson's class of a
hypersurface should be computable by an operation as simple as the
above. The original definition of MacPherson's class of $X$ involves
adding Mather's classes of subvarieties of $X$, weighted according to
their own singularity as measured by `local Euler obstructions'; the
one@-step (***) bypasses these complications.

We do not know whether this is an essential feature of hypersurfaces,
or whether a formula similar to (***) may compute
MacPherson's class of arbitrary varieties. While this is a natural
question, the approach of this paper does not seem well suited to
address it. One reason for presenting several different formulations of the
result is indeed that we do not see at the moment which formulation is
more suited to be generalized to arbitrary varieties.

These statements are given in \S1, Theorem I.1--5, where we offer several
alternative formulations for the class $c_*(X)$. More specifically:

---an explicit form in terms of Segre classes is given in \S1.1; 

---in \S1.2 we give a form which emphasizes the relation of MacPherson's
class with Fulton's; this was conjectured in \cite{Aluffi2}, where we
were putting restrictions on the line bundle $\Cal L=\Cal O(X)$ of $X$
and we could only prove the conjectured equality `numerically', that is
after taking degrees with respect to $\Cal L$. The result proved in this
note holds in the Chow group of $X$, and without hypotheses on $\Cal
O(X)$;

---in \S1.3 we discuss the form reproduced above; 

---in \S1.4 we introduce notations which we use elsewhere and in the
proof of the main Theorem, and give a formulation of the statement in
term of these notations;

---finally, in \S1.5 we express $c_*(X)$ as an explicit formula in terms of
a $\mu$@-class we introduced in \cite{Aluffi1}. In this form, the result is:
$$\cmp(X)=c(TM)\cap s(X,M)+c(\Cal L)^{\dim X}\cap (\mu_{\Cal
L}(Y)^\vee\otimes_M \Cal L)$$
(this formula also uses the notations introduced in \S1.4). In fact 
this is the form in which we first produced the result; the key technical
step in proving the main result of this paper can be written as a formula
detailing the behavior of the $\mu$@-class under blow@-ups (see \S 4.2).

The equivalence of the different formulations is essentially elementary,
and is proved along the way in \S1. The proof of the main result is in
\S\S2 and 3, and is on the contrary more delicate than this writer
originally expected. We show that a class satisfying the functorial
properties prescribed by Grothendieck and Deligne must necessarily agree
with $c_*(X)$; MacPherson's class $\cmp(X)$ is such a class. Our proof uses
resolution of singularities, reducing the statement to a computation for a
divisor with normal crossing (\S 2), and to studying the behavior of the
class under blow@-ups along nonsingular centers (\S 3). This ultimately
involves a careful analysis of an adaptation of MacPherson's graph
construction applied to our situation: the required equality is expressed
in terms of the vanishing of the contribution of a `cycle at infinity'
obtained in the construction, and this contribution can be evaluated
explicitly.

\S4 collects a few consequences of the main theorem. For example, in
its raw form (presented in \S1.1), the result writes $\cmp(X)$ as a
certain combination of Segre classes of the singular subscheme of $X$
in the ambient variety $M$. The good functoriality properties of $\cmp(X)$
translate then into properties of these Segre classes, which often we
are not able to prove more directly. We give one example of such
properties in \S~4.3, although this is bound to appear rather
technical to a reader who is not particularly fond of Segre classes.

However, properties of Segre classes often have very concrete geometric
applications. Here is an example of such an application, which can be
proved easily as a consequence of the main Theorem (see \S 4.4 for
details):
\proclaim{Proposition} Let $M_1$, $M_2$ be two hypersurfaces in $\P^N$, and
assume the contact scheme of $M_1$ and $M_2$ is smooth and
positive@-dimensional. Then $\deg M_1=\deg M_2$.\endproclaim 
Here the `contact scheme' is the singular subscheme of $M_1\cap M_2$ (in
the above sense), a scheme supported on the locus where $M_1$ and $M_2$ are
tangent. This corollary thus states that if two hypersurfaces of projective
space are {\it tangent\/}---in a strong sense---along a nonsingular
positive@-dimensional locus, then their degree must be equal: to our
knowledge, this observation is new.\vskip 6pt

The paper is organized so that a hasty reader who is willing to trust us on
the main technical step of the proof can skip \S3 at first: the notations
introduced in \S3 are not used elsewhere. In fact, an even more trusting
reader could skip all but the beginning of \S2.

Notations are gravely abused in this paper (pull@-backs are usually
omitted, etc.). The worst abuse occurs when we ask the reader to interpret
a class defined a priori on an ambient variety $M$, but supported on a
subvariety $X$, as a class defined on $X$: this will mean that there is
only one reasonable way to interpret the class as a push@-forward of a
class from $A_*X$, and that the given formula is a short@-hand for the
latter class. This saves us a great amount of notational grief, especially
in \S3; the reader who feels uncomfortable about this choice will only
believe our result after push@-forward to the ambient variety $M$.

I thank Tatsuo Suwa for useful comments.


\head \S1. Statements of the result\endhead

\subhead \S1.1. MacPherson's Chern class and Segre classes\endsubhead
A {\it constructible function\/} on an algebraic variety $X$ is an
integral linear combination of characteristic functions of closed
subvarieties. Over the complex numbers one can define a functor
assigning to $X$ the group of constructible functions on $X$; for any
proper morphism $f:X @>>> Y$, the push@-forward $f_*$ is defined by
setting
$$f_*(1_V)(p)=\chi(f^{-1}(p)\cap V)$$
for $V$ a subvariety of $X$ and $p\in Y$, and extending by linearity;
here $\chi$ denotes topological Euler characteristic. Grothendieck and
Deligne conjectured, and MacPherson proved (\cite{MacPherson}), that
there exists a natural transformation from this functor to homology,
such that, for $X$ nonsingular, $1_X$ is sent to the total homology
Chern class of $X$
$$c(TX)\cap [X]$$
It is then natural to consider the image of $1_X$ for arbitrary $X$; we
call this class {\it MacPherson's Chern class} of $X$, denoted $\cmp(X)$.
Some authors prefer the term {\it Schwartz@-MacPherson class\/} of
$X$, since the class was later shown to agree with a class previously
defined by M.--H.~Schwartz (\cite{Schwartz, B-S}). After work of
C.~Sabbah and G.~Kennedy (\cite{Kennedy}), this definition can be
extended to varieties over arbitrary fields of characteristic 0; and
$\cmp(X)$ lives in the Chow group $A_*X$ of cycles modulo rational
equivalence. Also, for possibly nonreduced $X$, we define
$\cmp(X)=\cmp(X_{\text{red}})$.

The main Theorem in this note will give explicit formulas for $\cmp(X)$,
in the case where $X$ is a hypersurface of an $n$@-dimensional
nonsingular variety $M$, that is, the zero@-scheme of a section $F$ of a
line bundle $\Cal L$ on $M$. We give the result in its raw form in terms of
Segre classes in this \S1.1; however, the result is more significant if the
formula is rewritten to emphasize its relation with other definitions of
characteristic classes for singular hypersurfaces. This is done in the rest
of this~\S1. In this section (\S 1.4) we will also introduce a notation
which makes some of the formulas easier to handle. The proof of the main
theorem will occupy \S\S 2 and 3.

Assume $\cmp$ satisfies the above functoriality condition, and assume
resolution of singularities \'a la Hironaka holds (for example, $\cmp$
can be MacPherson's class, in characteristic zero). By the {\it
singular (sub)scheme} of a hypersurface $X\subset M$ as above we mean
the subscheme of $M$ defined locally by $F$ and its first partial
derivatives---that is, by the {\it jacobian ideal} of $X$. (Note: the
definition of singular scheme of $X$ given in \cite{Aluffi1}, \S1.1,
fails to include $F$ among the local generators of its ideal. The
definition given here is the `correct' one; all the results in
\cite{Aluffi1} hold for this notion.) The singular scheme of $X$
depends on $X$ only, and not on a specific realization of $X$ as a
hypersurface of a nonsingular variety.

For $Y\subset X$, we define a class $s(X\setminus Y,M)\in A_*X$ by
setting its dimension@-$m$ component to be
$$s(X\setminus Y,M)_m=s(X,M)_m+(-1)^{n-m}\sum_{j=0}^{n-m}\binom {n-m}j
X^j\cdot s(Y,M)_{m+j}$$
(here and in the following $s(Y,M)$ denotes the {\it Segre class} of $Y$
in $M$ in the sense of \cite{Fulton}, Chapter 4).

\proclaim{Theorem I.1} Let $X$ be a hypersurface in a nonsingular
variety $M$, and let $Y$ be its singular scheme. Then $\cmp(X)=c_*(X)$,
where
$$c_*(X)=c(TM)\cap s(X\setminus Y,M)$$
\endproclaim

The expert will notice a similarity among the class $s(X\setminus
Y,M)$ defined above and formulas for residual intersections
(cf.~\cite{Fulton}, \S 9.2). Of course this is not accidental; the
connection will be clarified in the next subsection.

\subhead \S1.2. MacPherson's Chern class and Fulton's Chern
class\endsubhead
Let $X$ be a scheme embeddable in a nonsingular variety $M$. Fulton
shows (\cite{Fulton}, 4.2.6) that the class
$$c(TM)\cap s(X,M)$$
only depends on $X$ (and not on the choice of embedding). We will call
this class {\it Fulton's Chern class} of $X$, denoted $c_F(X)$.

Letting now $X$ again denote a divisor in a nonsingular variety $M$ of
dimension $n$, $\Cal L$ its line bundle, and letting $Y$ be its singular
subscheme, fix an integer $k\ge 0$ and consider the scheme $X(k)$
defined by `thickening $X$ $k$ times along $Y$': more precisely, if
$\Cal I_Y$ denotes the ideal of $Y$ and $\Cal J$ is the locally
principal ideal of $X$, then $X(k)$ is the subscheme of $M$ defined by
the ideal $\Cal J\cdot \Cal I_Y^k$. We may then consider the class
$$c_F(X(k))$$
in $A_*X$. We observed in \cite{Aluffi2} that this class is a {\it
polynomial} in $k$ with coefficients in $A_*X$, so it can be formally
evaluated for arbitrary $k$. It is also clear from the definition
that
$$c_F(X)=c_F(X(0))\quad.$$
Under the same hypotheses of the first statement, the main Theorem of
this note can then be stated as:
\proclaim{Theorem I.2} Let $X$ be a hypersurface in a nonsingular
variety $M$, and let $c_F(X(k))$ be defined as above. Then
$\cmp(X)=c_*(X)$, where
$$c_*(X)=c_F(X(-1))$$
\endproclaim

To see that this statement is equivalent to the first one given above
amounts to applying standard residual intersection formulas to compute
$s(X(k),M)$: by Proposition~9.2 in \cite{Fulton}, the $m$@-dimensional
component of this class is
$$s(X(k),M)_m=s(X,M)_m+\sum_{j=0}^{n-m}\binom {n-m}j
(-X)^j\cdot k^{n-m-j} s(Y,M)_{m+j}$$
That is, $s(X(-1),M)$ is the class $s(X\setminus Y,M)$ introduced
above, and it follows that the two expressions for $c_*(X)$ in the two
statements of the Theorem agree.

\subhead \S1.3. MacPherson's Chern class and Mather's Chern
class\endsubhead
Another notion of Chern classes for possibly singular varieties can be
defined as follows. For a reduced pure@-dimensional $X$ embedded in a
nonsingular $M$, let $X^\circ$ denote the nonsingular part of $X$. The
{\it Nash blow@-up} $\NB$ of $X$ is the closure in $\Grass_{\dim
X}(TM|_X)$ of the image of the map associating with every $p\in X^\circ$
the tangent space to $X^\circ$ at $p$; it comes equipped with a natural map
$\pi$ to $X$. The universal subbundle $\Cal T$ of $\Grass_{\dim
X}(TM|_X)$ restricts to $TX^\circ$ over $X^\circ\subset \NB(X)$, so it is
rather natural to consider the class
$$\pi_*(c(\Cal T)\cap [X])$$
Again, one proves this is independent of the choice of $M$; this class
is called {\it Mather's Chern class} of $X$, and we will denote it
$c_{M}(X)$. MacPherson computes $\cmp(X)$ as a suitable combination of
Mather's Chern classes of subvarieties of $X$.

The definition of $c_M(X)$ can be rewritten in a slightly different
way. If $X$ is a hypersurface in $M$, $x_1,\dots,x_n$ are local
parameters for $M$, and $F=0$ is a local equation for $X$, then
$\Grass_{\dim X}(TM)=\P T^*M$, the projectivized cotangent bundle to
$M$, and the map from $X^\circ$ is written in natural coordinates 
$$p \mapsto \left(p;\ddfi 1|_p:\dots:\ddfi n|_p\right)$$ 
It follows that locally $\NB$ can be recovered as the blow@-up of $X$
along the ideal of partials of the section of $\Cal L=\Cal O(X)$
defining $X$ in $M$ (cf.~also \cite{Nobile}, where a corresponding
observation is worked out for general $X$). We need to `globalize'
this observation.

\proclaim{Lemma I.1} Let $X$ be a hypersurface of a nonsingular variety
$M$, and let $Y$ be the singular subscheme of $X$. Then $\NB$ is the
proper transform of $X$ in the blow@-up $B\ell_YM @>\pi>> M$ of $M$
along $Y$. Further, there is a bundle $\Cal T$ on $B\ell_YM$ such that
$$c_M(X)=\pi_*(c(\Cal T)\cap [\NB])$$
\endproclaim
\demo{Proof} Let $\bdll$ denote the {\it bundle of principal parts of
$\Cal L$ over $M$,\/} where $\Cal L=\Cal O(X)$ as above, and consider
the section $M @>>> \bdll$ determined by $F\in H^0(M,\Cal L)$; locally
we can write this as
$$p \mapsto \left(p; F(p): \ddfi 1|_p:\dots:\ddfi n|_p\right)$$
(recall that $\bdll$ fits in an exact sequence
$$ 0 @>>> T^*M\otimes \Cal L @>>> \bdll @>>> \Cal L @>>> 0\quad.)$$
It follows that the closure of the image of the corresponding rational
map
$$M \dashrightarrow \P\bdll$$
is the blow@-up $B\ell_YM @>\pi >> M$ along the singular scheme $Y$ of
$X$ defined in \S1.1. Over $X$ this reproduces (up to tensoring by $\Cal
L$) the map to $\P T^*M$ considered above, and the first part of the
statement follows.

For the second part, let $\Cal Q$ be the universal quotient bundle
of $\P\bdll$, and consider $\Cal T=\Cal Q^\vee\otimes\Cal L$. Observe
that the universal subbundle $\Cal O(-1)$ of $\P\bdll$ is $\subset
T^*M\otimes \Cal L$ over $\NB$ (since this dominates $X$, and $F\equiv 0$
over $X$), so there is a sequence
$$ 0 @>>> \Cal Q'|_{\NB} @>>> \Cal Q|_{\NB} @>>> \Cal L @>>> 0 $$
where $\Cal Q'$ is the universal quotient bundle of $\P (T^*M\otimes
\Cal L|_X)$. Chasing the identification 
$$\Grass_{n-1}(TM|_X)\cong \P(T^*M|_X) \cong \P(T^*M\otimes \Cal L|_X)$$
shows that ${\Cal Q'}^\vee \otimes\Cal L$ is the bundle used above in
the definition of $c_M$, and $\Cal T=\Cal Q^\vee \otimes\Cal L$ differs
from this by a trivial factor.\qed\enddemo

The next Lemma gives Mather's class in the form stated in the
introduction. In $B\ell_YM$ we have the
proper transform $\NB$ of $X$; $X$ itself pulls back to a Cartier divisor
of $B\ell_YM$, which we still denote $X$; also, we have the exceptional
divisor $\Cal Y$, and we note that $\Cal O(\Cal Y)$ is the restriction
of $\Cal O(-1)$ from $\P\bdll$.
\proclaim{Lemma I.2} With notations as above,
$$c_M(X)=c(TM)\cap \pi_*\left(\frac{[\NB]}{1+X-\Cal Y}\right)$$
\endproclaim
\demo{Proof} From the proof of Lemma~I.1, we know
$$c_M(X)=\pi_*\left(c\left(\left(\frac{\bdll}{\Cal
O(-1)}\right)^\vee\otimes\Cal L\right) \cap [\NB]\right)$$
Using $c(\Cal L)=1+X$, $c(\Cal O(1))=1-\Cal Y$, and
$c(\bdll)=c(T^*M\otimes \Cal L)c(\Cal L)$, this gives
$$\align
c_M(X) &=\pi_*\left(\frac{c((\bdll)^\vee \otimes\Cal L)}{1+X-\Cal Y}
\cap [\NB]\right)\\
&=\pi_*\left(c(TM\otimes\Cal L^\vee \otimes \Cal L)\,c(\Cal
L^\vee\otimes\Cal L)\cap \frac{[\NB]}{1+X-\Cal Y}\right)\\
&=c(TM)\cap \pi_*\left(\frac{[\NB]}{1+X-\Cal Y}\right)\quad\qed
\endalign$$
\enddemo

With this understood, the third form of the result of this paper should
appear more significant; this is the form stated in the introduction.
Under the same hypotheses of the first and second statement of the main
Theorem:
\proclaim{Theorem I.3} Let $X$ be a hypersurface in a nonsingular
variety $M$, and let $\pi$, $\Cal Y$ be as defined above. Then
$\cmp(X)=c_*(X)$, where
$$c_*(X)=c(TM)\cap \pi_*\left(\frac{[X]-[\Cal Y]}{1+X-\Cal Y}\right)$$
\endproclaim

To see that this is equivalent to the other statements, we have to show
that
$$\pi_*\left(\frac{[X]-[\Cal Y]}{1+X-\Cal Y}\right)=s(X\setminus Y,M)
\tag*$$
One way to see this is to consider for all $k$ the class
$$\pi_*\left(\frac{[X]+k\,[\Cal Y]}{1+X+k\,\Cal Y}\right)\quad:$$
for $k\ge 0$, $[X]+k\,[\Cal Y]$ is the cycle of the inverse image of
the scheme $X(k)$ considered in \S1.2, hence by the birational
invariance of Segre classes
$$\pi_*\left(\frac{[X]+k\,[\Cal Y]}{1+X+k\,\Cal Y}\right) = s(X(k),M)$$
for $k\ge 0$. Both sides are polynomials in $k$, so they must agree for
$k=-1$, and this is what (*) claims.

\subhead \S1.4. Notational device\endsubhead
The following notations will be of help in writing the arguments needed
to prove the Theorem stated above. In fact they simplify considerably
the first formulation we gave, by giving a summation@-free alternative
definition for the class $s(X\setminus Y,M)$. 

\definition{Definition} If $A=\oplus_i a^i$ is a rational equivalence
class on a scheme, indexed by codimension, we let
$$A^\vee=\sum_{i\ge 0} (-1)^i a^i\quad,$$
the {\it dual} of $A$; also, for a line bundle $\Cal L$ we let
$$A\otimes \Cal L=\sum_{i\ge 0} \frac{a^i}{c(\Cal L)^i}\quad,$$
the {\it tensor of $A$ by $\Cal L$.} We put a subscript to
$\otimes$ to denote the ambient in which the tensor is taken, if this
doesn't seem otherwise clear from the context.\enddefinition

These notations were introduced in \cite{Aluffi2}, where we also proved
simple compatibilities with standard vector bundle operations (Prop.~1
and 2 in \cite{Aluffi2}, \S2). We will freely use those properties in
this note (especially in \S2).

\proclaim{Lemma I.3} With notations as above,
$$s(X\setminus Y,M)=s(X,M)+c(\Cal L)^{-1}\cap(s(Y,M)^\vee \otimes \Cal
L)$$
\endproclaim
\demo{Proof} This is detailed in \cite{Aluffi2} (section 2), so we will
not reproduce it here. It is a good exercise for the reader
interested in acquiring some familiarity with the notations introduced
above: if $X\subset W\subset M$, with $X$, $M$ as above, and $R$ is the
residual scheme to $X$ in $W$, show that
$$s(W,M)=s(X,M)+c(\Cal L)^{-1}\cap (s(R,M)\otimes\Cal L)$$
Applying to $X(k)$ (as in \S 1.2 above) and setting $k=-1$ gives the
statement. Details may be found in \cite{Aluffi2}.\qed\enddemo

By Lemma I.3, the main Theorem can be formulated:
\proclaim{Theorem I.4} Let $X$ be a hypersurface in a nonsingular
variety $M$, and let $Y$ be its singular scheme. Then $\cmp(X)=c_*(X)$,
where
$$c_*(X)=c(TM)\cap\left(s(X,M)+c(\Cal L)^{-1}\cap(s(Y,M)^\vee \otimes_M
\Cal L)\right)$$
\endproclaim

\subhead \S1.5. MacPherson's Chern class and $\mu$@-classes\endsubhead
Another formulation of the main result of this paper can be given in
terms of the $\mu$@-class introduced in \cite{Aluffi1}.

If $Y$ is a singular subscheme of a hypersurface $X$ in a nonsingular
variety $M$, with $\Cal L=\Cal O(X)$, the {\it $\mu$@-class of $Y$ with
respect to $\Cal L$} is defined by
$$\mu_{\Cal L}(Y)=c(T^* M\otimes\Cal L)\cap s(Y,M)$$
This class does not depend on the specific realization of $Y$ as a
singular scheme of a hypersurface (Corollary 1.7 in \cite{Aluffi1}).
Note that if $Y\subset M$ is nonsingular, and with the notations
introduced in \S1.4:
$$c(\Cal L)^{\dim M}\cap \left(\mu_{\Cal L}(Y)^\vee\otimes_M \Cal
L\right)=(-1)^{\codim_M Y}c(TY)\cap [Y]$$
(this isn't entirely obvious from the definition; it follows from
Corollary 1.8 in \cite{Aluffi1}. Of course this formula fails
spectacularly unless one assumes $Y$ is realized as the singular scheme
of a hypersurface in $M$.) That is, the $\mu$@-class can be used to
define yet another class extending the notion of Chern class to (certain)
possibly singular varieties. It turns out that this class gives a
precise `correction term' for the class of the hypersurface of which
$Y$ is the singular scheme. The precise statement is
\proclaim{Theorem I.5} Let $X$ be a hypersurface in a nonsingular
variety $M$, let $Y$ be its singular scheme, and $\Cal L=\Cal
O(X)$. Then $\cmp(X)=c_*(X)$, where
$$c_*(X)=c(TM)\cap s(X,M)+c(\Cal L)^{\dim X}\cap (\mu_{\Cal
L}(Y)^\vee\otimes_M \Cal L)$$ 
\endproclaim
This statement is equivalent to Theorem~I.1--4; showing this amounts to
showing that
$$c(\Cal L)^{\dim X}\cap\left(\mu_{\Cal L}(Y)\otimes_M \Cal L\right) =
c(TM)\cap s(X\setminus Y,M)\quad.$$
This is a good exercise in the notations of \S1.4, and we leave it to
the reader.\vskip 6pt

The proof of Theorem I.1--5 occupies the next two sections. \S2 reduces
it to showing that the class $c_*$ introduced in this section satisfies
a simple blow@-up formula ((3) in \S2); \S3 proves this
formula.


\head \S2. The proof: preliminaries\endhead
In this section we set up our strategy for the proof of the main
Theorem, by reducing it to the proof of one property (property (3))
below describing the behavior of $c_*$ under blow@-ups along
nonsingular subvarieties. Section 3 is devoted to the proof of
(3).

We will need here neither the full functorial picture summarized in the
beginning of \S1 nor the details of MacPherson's construction. We will
only consider the following three properties of a class $c_*$:\roster
\item If $X$ is a hypersurface, and the support $\red X$ of $X$ is
nonsingular, then

$$c_*(X)=c(T\red X)\cap [\red X]\quad;$$
\item {\it Normal crossings:\/} if $X=X_1\cup X_2$, where $X,X_1$
are divisors with normal crossings in a nonsingular ambient variety $M$,
and $X_2$ is a nonsingular hypersurface of $M$, then
$$c_*(X)=c_*(X_1)+c_*(X_2)-c_*(X_1\cap X_2)\quad;$$
\item {\it Blow@-up:\/} if $X$ is a hypersurface of a nonsingular
variety $M$, $Z\subset X\subset M$ is a nonsingular subvariety of
codimension $d$ in $M$, $\pi:\Til M @>>> M$ is the blow@-up of $M$
along $Z$, and $X'$ denotes the (scheme@-theoretic) inverse image of $X$
in $\Til M$, then
$$(\pi_{|X'})_*(c_*(X'))=c_*(X)+(d-1)\,c_*(Z)\quad.$$
\endroster

\noindent (We often abuse notations and omit obvious push@-forwards and
pull@-backs, as above.) It is easy to see that MacPherson's class
satisfies these properties: (2) follows from the analogous relation
between characteristic functions (in fact, with no restrictions on
what $X_1$, $X_2$ may be); for (3), consider the map
$$f:X'\amalg Z @>>> X$$
restricting to $\pi_{|X'}$ on $X'$ and to the inclusion $i$ into $X$ on
$Z$; and define the constructible function $\aleph=1_{X'}-(d-1)\,1_Z$ on
$X'\amalg Z$. Then for $p\in X$
$$f_*(\aleph)(p)=(\pi_{|X'})_*(1_{X'})(p)-(d-1)\,i_*(1_Z)(p)
=\left\{\aligned 1-0=1\quad &\text{for $p\notin Z$}\\
d-(d-1)=1\quad &\text{for $p\in Z$}
\endaligned\right.$$
since for $p\in Z$, $\pi^{-1}(p)\cong \P^{d-1}$ has Euler characteristic
$d$. This shows $f_*(\aleph)=1_X$, and (3) follows for $\cmp$.

It is also clear from embedded resolution of singularities that the
class is uniquely determined by (1), (2), (3) for hypersurfaces $X$ of
nonsingular varieties $M$. For this, let $M_s @>>> M_{s-1} @>>> \dots
@>>> M_1 @>>> M_0=M$ be a sequence of blow@-ups at non singular centers
such that the inverse image $X_s$ of $X_0=X$ is a divisor with
nonsingular components and normal crossings. Then (1) and (2) determine
the class for $X_s$; and at each stage $X_i @>>> X_{i-1}$ is a map as in
(3), so the value of the class at $X_{i-1}$ is determined by its value
at $X_i$, $i=1,\dots,s$.

Summarizing: {\it in order to prove the main Theorem, it suffices to
show that the class $c_*(X)$ introduced in the statements of 
Theorem~I in \S1 satisfies properties (1), (2), (3) above.\/}

The rest of this section is devoted to the proof of (1) and (2) for
this class. Property (3) is technically more demanding, and we will
devote the entire \S3 to its proof.

\subhead \S2.1. $c_*(X)=c_*(\red X)$\endsubhead
Here we prove that under good hypotheses the hypersurface may be
assumed to be {\it reduced.\/} Note: it will be a consequence of the
main Theorem that in fact this holds for {\it all\/} hypersurfaces, but
we do not know how to prove this directly in general. The following 
lemma suffices for (1) above, and simplifies the work required to prove
(2).

The context of the lemma is as follows: we want to show that if the
components of a hypersurface are sufficiently transversal, then the
class of their union is independent of the multiplicity with which the
components appear.
\proclaim{Lemma II.1} Let $X_1$, $X_2$ be two hypersurfaces in a
nonsingular variety $M$, and assume that at every point of $M$ there are
local parameters $x_1,\dots,x_n$ such that $X_1$ has equation $x_1=0$,
and $X_2$ has equation $f(x_2,\dots,x_n)=0$. Also denote by $X_{(m)}$
the hypersurface locally defined by the ideal $(x_1^m f)$. Then
$$c_*(X_{(m)})=c_*(X_{(1)})\quad\text{for all $m\ge 1$}$$
\endproclaim
\demo{Proof} The jacobian ideal of $X_{(m)}$ is given locally by
$$x_1^{m-1}\left(f,x_1\ddfi i\right)_{i\ge 2}\quad:$$
that is, it consists of the $(m-1)$@-multiple of $X_1$ and of a
residual $R$ independent of $m$. But then
$$B\ell_{Y_{(m)}} M\cong B\ell_R M @>\pi>> M$$
is independent of $m$, while the exceptional divisor
in $B\ell_{Y_{(m)}}$ is 
$$\Cal Y_{(m)}=(m-1) X_1+ \Cal Y$$
where $\Cal Y$ denotes the exceptional divisor in $B\ell_R M$. Using
the expression for $c_*$ in Theorem~I.3:
$$\align
c_*(X_{(m)})&=c(TM)\cap \pi_*\frac{[X_{(m)}]-[\Cal Y_{(m)}]}
{1+X_{(m)}-\Cal Y_{(m)}}\\
&=c(TM)\cap \pi_*\frac{(m [X_1]+[X_2])-((m-1)[X_1]+[\Cal Y])}
{1+(m X_1+X_2)-((m-1)X_1+\Cal Y)}\\
&=c(TM)\cap \pi_*\frac{[X_1]+[X_2]-[\Cal Y]}{1+X_1+X_2-\Cal Y}
\endalign$$
is also independent of $m\ge 1$, as needed.\qed\enddemo
Lemma II.1 implies (1):
\proclaim{Corollary II.1} If the support $\red X$ of $X$ is nonsingular,
then
$$c_*(X)=c(T\red X)\cap [\red X]$$
\endproclaim
\demo{Proof} Taking $X_2=\emptyset$ in Lemma II.1 yields
$c_*(X)=c_*(\red X)$: so we may assume $X$ is reduced and nonsingular.
Then its singular scheme is $Y=\emptyset$, so $s(X\setminus
Y,M)=s(X,M)$. Finally $s(X,M)=c(N_XM)^{-1}\cap [X]$ (the inverse Chern
class of the normal bundle of $X$), so
$$c_*(X)=c(TM) c(N_XM)^{-1}\cap [X]=c(TX)\cap [X]\quad\qed$$
\enddemo

By a {\it divisor with normal crossings} we mean a union of smooth
distinct hypersurfaces $X_1\cup\dots\cup X_r$ such that at each point
of intersection of some of the $X_i$, say of $X_1,\dots,X_k$, there are
local coordinates $(x_1,\dots,x_n)$ for the ambient variety so that
$x_1=0,\dots,x_k=0$ are equations for $X_1,\dots,X_k$ respectively. In
fact we must allow the $X_i$'s to come with multiplicity: the plan is
to apply resolution of singularities to an arbitrary hypersurface
$X\subset M$, and this will produce a nonsingular variety mapping to
$M$, in which the (scheme@-theoretic) inverse image of $X$ is a
divisor with normal crossings, whose components will appear with
multiplicity. Lemma II.1 implies that at this stage we will be able to
discard the extra multiplicity information:
\proclaim{Corollary II.2} If $X$ is a (possibly nonreduced) divisor with
normal crossing, then
$$c_*(X)=c_*(\red X)$$
\endproclaim
\demo{Proof} This follows by repeatedly applying Lemma II.1, taking for
$X_1$ each component of the divisor in turn.\qed\enddemo

It will follow from the main Theorem that in fact Corollary II.2
holds for arbitrary hypersurfaces. Again, we do not know how to prove
directly this more general statement.

\subhead \S2.2. Divisors with normal crossings: proof\endsubhead
Here we prove that $c_*$ satisfies property (2) above. We first
translate (2) into the exact form proved below.

Let $X=X_1\cup\dots\cup X_r$ be a divisor with normal crossings. By
Corollary~2 above, in computing $c_*(X)$ we may assume $X$ is reduced.
As usual, $Y$ denotes the singular scheme of $X$ and $\Cal L=\Cal
O(X)$. Also, we write $\Cal L_i$ for $\Cal O(X_i)$.

\proclaim{Lemma II.2} In order to prove (2), it suffices to show that
$$s(Y,M)=\left(\left(1-\frac{c(\Cal L^\vee)}{c(\Cal L_1^\vee)\cdots
c(\Cal L_r^\vee)}\right)\cap [M]\right)\otimes \Cal L$$
\endproclaim
\noindent Here the reader must interpret the right@-hand@-side as a class
supported in $Y$, that is, obvious cancellations must be performed on the
right@-hand@-side. This will be assumed implicitly in the following.

\demo{Proof} Assuming $s(Y,M)$ is given by the expression in the statement,
we derive
$$s(Y,M)^\vee\otimes\Cal L =\left(1-\frac{c(\Cal L)}{c(\Cal L_1)\cdots
c(\Cal L_r)}\right)\cap [M]$$
and therefore (using the expression for $c_*$ given in Theorem~I.4, and
after simple manipulations)
$$c_*(X)=c(TM)\cdot \left(1-\frac 1{(1+X_1)\cdots (1+X_r)}\right)\cap
[M]$$
Thus showing (2) for this class amounts to showing that
$$\multline
c(TM)\left(1-\frac 1{(1+X_1)\cdots (1+X_r)}\right)\cap [M] =
c(TX_1)\cap [X_1] +c(TM)\cdot\\
\cdot\left(1-\frac 1{(1+X_2)\cdots (1+X_r)}\right)\cap [M]
-c(TX_1)\left(1-\frac 1{(1+X_2)\cdots (1+X_r)}\right)\cap [X_1]
\endmultline$$
since $X_1\cap (X_2\cup \cdots\cup X_r)$ is also a reduced divisor with
normal crossings (in $X_1$, which is assumed to be nonsingular). Now
the right@-hand@-side can be written
$$c(TM)\left(\frac{X_1}{(1+X_1)\cdots (1+X_r)}+1- \frac
1{(1+X_2)\cdots (1+X_r)}\right)\cap [M]$$
and this is immediately seen to equal the left@-hand@-side, as
needed.\qed\enddemo

In passing we note that since $\cmp$ satisfies (2), the following
formula must hold for MacPherson's Chern class of a reduced divisor
$X=X_1\cup\cdots\cup X_r\subset M$ as above:
$$\cmp(X)=c(TM)\cdot \left(1-\frac 1{(1+X_1)\cdots (1+X_r)}\right)\cap
[M]\quad.$$

By Lemma II.2, we are reduced to showing
$$s(Y,M)=\left(\left(1-\frac{c(\Cal L^\vee)}{c(\Cal L_1^\vee)\cdots
c(\Cal L_r^\vee)}\right)\cap [M]\right)\otimes \Cal L$$
for $Y$ the singular scheme of a reduced divisor $X=X_1\cup\cdots\cup
X_r$ with nonsingular components and normal crossings.

\demo{Proof of (2)} As a set, $Y$ is the union of all the $X_i\cap X_j$ with
$i\ne j$; $X$ has multiplicity $k$ along the intersection of $k$
components, $X_I=X_{i_1}\cap \dots\cap X_{i_k}$ ($k=|I|$), provided that
this is nonempty. We will work by induction on the number $N$ of
nonempty such intersections $X_I$, $|I|\ge 2$.

The statement is clear if this number is 0, that is if $Y$ is empty:
$c(\Cal L_1^\vee\otimes\cdots\otimes \Cal L_r^\vee) = c(\Cal
L_1^\vee)\cdots c(\Cal L_r^\vee)$ if the $X_i$'s do not intersect.
Assume then $Y\ne\emptyset$, and consider an $X_I$ of minimal dimension,
say $Z=X_1\cap\dots\cap X_k$. Locally along $Z$, $X_i$ has then equation
$x_i=0$ (for $i\le k$), where the $x_i$'s are part of a system of
parameters; so the hypersurface is $x_1\cdots x_k=0$ along $Z$, and $Z$
has (local) ideal $(x_1,\dots,x_k)$; along $Z$, $Y$ has ideal
$$\left(\frac{x_1\cdots x_k}{x_1},\dots,\frac{x_1\cdots x_k}{x_k}
\right)$$ Note that if some other hypersurface of the lot came in at
some point of $Z$ not covered by the above chart, this would determine
a smaller nonempty intersection, against the minimality of $Z$. In
other words, $X_i\cap Z=\emptyset$ for $i>k$.

Now blow@-up $M$ along $Z$; with a suitable choice of coordinates
$\tilde x_i$ in the blow@-up, we can write the blow@-up map as
$$\left\{\aligned
x_1&=\tilde x_1\\
x_2&=\tilde x_1\tilde x_2\\
&\dots\\
x_k&=\tilde x_1\tilde x_k
\endaligned\right.$$
and the inverse image of $Y$ has ideal
$$\left(\frac{\tilde x_1^k\cdots \tilde x_k}{\tilde x_1},\frac{\tilde
x_1^k\cdots \tilde x_k}{\tilde x_1\tilde x_2},\dots,\frac{\tilde
x_1^k \cdots\tilde x_k}{\tilde x_1\tilde x_k}\right) = \tilde x_1^{k-1}
\left(\frac{\tilde x_2\cdots \tilde x_k}{\tilde x_2},\dots,\frac{\tilde
x_2\cdots\tilde x_k}{\tilde x_k}\right)$$
in this chart. Now this says that the residual of $(k-1)$ times the
exceptional divisor in the inverse image of $Y$ is (in this chart) the
{\it singular scheme of the proper transform of the hypersurface.\/}
This must in fact hold globally on $Z$, as the behavior on the other
charts is identical to the one shown above.

Now the key is that the proper transform of the hypersurface is again
a divisor with normal crossing, but for which the number $N$
considered above is one less than for the original hypersurface;
therefore by induction we know the Segre class of its singular scheme:
{\eightpoint
$$1-\frac{1-(X_1-E)-\dots-(X_k-E)-X_{k+1}-\dots-X_r}
{(1-X_1+E)\dots(1-X_k+E)(1-X_{k+1})\dots(1-X_r)}\otimes \Cal
O(X_1+\dots+X_r-kE)$$}
\noindent where $E$ is the class of the exceptional divisor.

Using Proposition 3 from \cite{Aluffi2} to throw in $(k-1)E$, and
using the birational invariance of Segre classes, we get that $s(Y,M)$
is the push@-forward to $M$ of{\eightpoint
$$\multline
\frac{(k-1)E}{(1+(k-1)E)}+\frac 1{1+(k-1)E}\cdot\\
\cdot\left(1-\frac{1-X_1-\dots-X_r+k\, E}
{(1-X_1+E)\dots(1-X_k+E)(1-X_{k+1})\dots(1-X_r)}\otimes \Cal
O(X_1+\dots+X_r-E)\right)
\endmultline$$}
\noindent that is{\eightpoint
$$\multline
1-\frac 1{1+(k-1)E}\cdot\\
\cdot\left(\frac{1-X_1-\dots-X_r+k\, E}
{(1-X_1+E)\dots(1-X_k+E)(1-X_{k+1})\dots(1-X_r)}\otimes \Cal
O(X_1+\dots+X_r-E)\right)
\endmultline$$}
and using \cite{Aluffi2}, \S2, this is manipulated into
{\eightpoint
$$\align
1- &\left(\frac{1-X_1-\dots-X_r+ E}{1-X_1-\dots-X_r+k E}\cdot\right.\\
&\left. \cdot\frac{1-X_1-\dots-X_r+k E}
{(1-X_1+E)\dots(1-X_k+E)(1-X_{k+1})\dots(1-X_r)}\otimes \Cal
O(X_1+\dots+X_r-E)\right)\\
=1- &\left(\frac{1-X_1-\dots-X_r+ E}
{(1-X_1+E)\dots(1-X_k+E)(1-X_{k+1})\dots(1-X_r)}\otimes \Cal
O(X_1+\dots+X_r-E)\right)\\
=1- &\left(\frac{(1-X_1-\dots-X_r)}{(1-X_1)\dots(1-X_k)}\cdot
\frac{(1-E)^{r-1}}{(1-X_{k+1}-E)\dots(1-X_r-E)}\otimes \Cal
O(X_1+\dots+X_r)\right)
\endalign
$$}
Now we claim that
$$\frac{(1-E)^{r-1}}{(1-X_{k+1}-E)\dots(1-X_r-E)}\tag*$$
pushes forward to
$$\frac 1{(1-X_{k+1})\dots(1-X_r)}\tag**$$
Indeed, any term involving both $E$ and some of the $X_i$'s, $i> k$,
is necessarily 0 since these $X_i$'s do not meet $Z$; so (*) equals
$$(1-E)^{k-1}-1+\frac 1{(1-X_{k+1})\dots(1-X_r)}\quad;$$
and all powers $E^i$ with $0<i <k$ push forward to 0 because $Z$ has
codimension $k$. So (**) is all that survives the push@-forward.

In conclusion, this shows that $s(Y,M)$ equals
$$[M]- \frac{([M]-[X_1]-\dots-[X_r])}{(1-X_1)\dots(1-X_r)} \otimes \Cal
O(X_1+\dots+X_r)\quad,$$
completing the induction step.\qed\enddemo\vskip 6pt

This concludes the proof that the class $c_*$ of \S\S1 and 2
satisfies properties (1) and (2) stated at the beginning of this
section (and it follows that $c_*$ and MacPherson's class coincide for
hypersurfaces with normal crossing).

{\it \S3 will be devoted to the proof of (3), thereby concluding the
proof of Theorem~I.\/} This will be by far the most delicate ingredient
in the proof of the main Theorem.


\head \S3. The proof: behavior under blow@-ups\endhead
The last ingredient in the proof of the main Theorem is the
proof of (3) from \S2. We will obtain this by transforming (3) into
equivalent and more basic assertions, which however will require more
and more notations to be stated. In the end, (3) will follow by an
explicit computation of a `cycle at infinity' ($Z_\infty$ in \S\S3.4--3.8)
arising in a graph construction.

We first reproduce the notations used so far, and the statement of (3)
given in \S2. Let $X$ be a hypersurface of a nonsingular variety $M$,
and let $Z\subset X\subset M$ be a nonsingular subvariety of codimension
$d$ in $M$. $\Til M @>\pi>> M$ will be the blow@-up of $M$ along $Z$, 
$E$ will denote the exceptional divisor of this blow@-up, and $X'$ the
scheme@-theoretic inverse image of $X$ in $\Til M$; $\Cal L$ will be the
line bundle of $X$ (hence its pull@-back, also denoted $\Cal L$, is the
line bundle of $X'$), and $Y$, $Y'$ will respectively denote the
singular schemes of $X$, $X'$. Then (3) states that
$$\pi_*(c_*(X'))=c_*(X)+(d-1)\,c_*(Z)$$
in $A_*X$. (Note: in this section especially we will often incur in
severe notational abuses, of which this formula is a good sample. To
interpret this formula correctly, the reader is expected to restrict
$\pi$ to $X'$ before using it to push@-forward $c_*(X')$; and to push
forward $c_*(Z)$ from $A_*Z$ to $A_*X$. While this will make some of our
statements slightly imprecise, employing full notations would often make
them quite unreadable; we opt for the first alternative.)

\subhead \S3.1. (3) in terms of classes in $\P\bdll$, $\P\bdlr$
\endsubhead
Here we translate (3) by using the form of $c_*$ given in Theorem I.3.

We will denote by $\bdll$, $\bdlr$ respectively the bundles of principal
parts of $\Cal L$ over $M$, $\Til M$. The section $F$ of $\Cal L$ over
$M$, $M'$ defining $X$, $X'$ resp.~determine sections
$$\Til M @>>> \pi^*\bdll\qquad;\qquad \Til M @>>> \bdlr$$
which projectivize to rational maps
$$\Til M \dashrightarrow \P\pi^*\bdll \qquad;\qquad \Til M
\dashrightarrow \bdlr$$
The closures of the images of these maps are the blow~ups
$B\ell_{\pi^{-1}Y} \Til M$, $B\ell_{Y'} \Til M$ respectively (this
follows from staring at local descriptions for the sections, cf.~\S1.3).

This is the first instance in which we perform two parallel
constructions: one on the $\pi^*\bdll$ side, the other on the $\bdlr$
side. As a rule, we will put subscripts $M$, $\Til M$ on corresponding
objects in the two sides, to keep track of which side they belong to: we
start this convention by naming the universal subbundles in
$\P\pi^*\bdll$, $\P\bdlr$ respectively $\Cal O_M(-1)$, $\Cal O_{\Til
M}(-1)$. Similarly, $\Cal Y_M$, $\Cal Y_{\Til M}$ will denote
respectively the exceptional divisors in $B\ell_{\pi^{-1}Y} \Til M$,
$B\ell_{Y'} \Til M$; note that $\Cal O(\Cal Y_M)$, $\Cal O(\Cal Y_{\Til
M})$ are respectively the restriction of $\Cal O_M(-1)$, $\Cal O_{\Til
M}(-1)$ to the blow@-ups. Also, $p_M$, $p_{\Til M}$ will denote
respectively the bundle maps on $\P\pi^*\bdll$, $\P\bdlr$. Finally,
the reader is warned that the $\pi^*$ employed so far will soon be
dropped (as is allowed by various functorialities of pull@-backs).

\proclaim{Claim III.1} In order to prove (3), it suffices to show that
$$\pi_*{p_M}_*\left(c\left(\frac{\bdll}{\Cal O_M(-1)}\right)\cap
[B\ell_{\pi^{-1}Y} \Til M]\right) = \pi_*{p_{\Til
M}}_*\left(c\left(\frac{\bdlr} {\Cal O_{\Til M}(-1)}\right)\cap
[B\ell_{Y'} \Til M]\right)$$
\endproclaim
\noindent(see the last paragraph in the introduction for
clarifications on the notations.)
\demo{Proof} Writing $c_*$ as in \S1.3, and with the above notations,
$$\gather
c_*(X)=\pi_*\left(c(TM)\cap {p_M}_*\left(\frac{[X]-[\Cal
Y_M]}{1+X-\Cal Y_M}\right)\right)\\
c_*(X')=c(T\Til M)\cap {p_{\Til M}}_*\left(\frac{[X']-[\Cal
Y_{\Til M}]}{1+X'-\Cal Y_{\Til M}}\right)
\endgather$$
Now{\eightpoint
$$\align
&\pi_*{p_M}_*\left(c\left(\frac{\bdll}{\Cal O_M(-1)}\right)\cap
[B\ell_{\pi^{-1}Y} \Til M]\right) - \pi_*{p_{\Til
M}}_*\left(c\left(\frac{\bdlr} {\Cal O_{\Til M}(-1)}\right)\cap
[B\ell_{Y'} \Til M]\right)\\
&=\pi_*\left({p_M}_*\left(\frac{c(T^*M\otimes\Cal L)c(\Cal L)}{1+\Cal
Y_M}\cap [B\ell_{\pi^{-1}Y} \Til M]\right) -
{p_{\Til M}}_*\left(\frac{c(T^*\Til M\otimes\Cal L)c(\Cal L)}{1+\Cal
Y_{\Til M}} \cap [B\ell_{Y'} \Til M]\right)\right)\\
&=c(\Cal L)^{n+1}\pi_*\left({p_M}_*\left(\frac{c(T^*M)}{1-X+\Cal
Y_M}\cap [B\ell_{\pi^{-1}Y} \Til M]\right) - {p_{\Til M}}_*
\left(\frac{c(T^*\Til M)}{1-X'+\Cal Y_{\Til M}} \cap
[B\ell_{Y'} \Til M]\right)\right)\otimes\Cal L
\endalign$$}
(notations as in \S1.4, and properties of the same from
\cite{Aluffi2}). Thus the equality in the statement is equivalent to
$$\pi_*\left({p_M}_*\left(\frac{c(T^*M)}{1-X+\Cal
Y_M}\cap [B\ell_{\pi^{-1}Y} \Til M]\right) - {p_{\Til M}}_*
\left(\frac{c(T^*\Til M)}{1-X'+\Cal Y_{\Til M}} \cap
[B\ell_{Y'} \Til M]\right)\right)=0$$
Taking duals, this is equivalent to
$$\pi_*\left({p_M}_*\left(\frac{c(TM)}{1+X-\Cal Y_M}\cap
[B\ell_{\pi^{-1}Y} \Til M]\right) - {p_{\Til M}}_*
\left(\frac{c(T\Til M)}{1+X'-\Cal Y_{\Til M}} \cap [B\ell_{Y'} \Til
M]\right)\right)=0$$ that is, to
$$\multline
\pi_*{p_M}_*\left(c(TM)\left(1-\frac{X-\Cal Y_M}{1+X-\Cal Y_M}\right)
\cap [B\ell_{\pi^{-1}Y} \Til M]\right)\\
- \pi_*{p_{\Til M}}_*\left(c(T\Til M)\left(1-\frac{X'-\Cal Y_{\Til
M}}{1+X'-\Cal Y_{\Til M}}\right) \cap [B\ell_{Y'} \Til M]\right)=0
\endmultline$$
or, using the expressions given above for $c_*$, to:{\eightpoint
$$
\pi_*{p_M}_*\left(c(TM)\cap [B\ell_{\pi^{-1}Y} \Til
M]\right)-c_*(X)
- \pi_*{p_{\Til M}}_* \left(c(T\Til M)\cap [B\ell_{Y'} \Til M]\right)
+\pi_*c_*(X')=0$$}
and finally, using the projection formula, to
$$\pi_*c_*(X')=c(X)+\pi_*\left((c(T\Til M)-c(TM))\cap [\Til M]\right)$$
Now
$$\pi_*(c(T\tilde M)\cap [\tilde M])-c(TM)\cap [M]=(d-1)\,c(TZ)\cap
[Z]\quad:$$
in characteristic 0 this is immediate from the functoriality of
MacPherson's Chern classes; but it holds in general, as may be easily
checked using Theorem 15.4 in \cite{Fulton}. Therefore the equality in
the statement is equivalent to
$$\pi_*c_*(X')=c(X)+(d-1)\,c(TZ)\cap [Z]$$
which is precisely (3), as needed.\qed\enddemo

\subhead \S3.2. (3) in terms of classes in $\P(\bdll\oplus\bdlr)$ over
$\Til M$\endsubhead
Before attacking the equality stated in Claim~III.1, we need another
notational layer to put both sides in the same place. The general plan
is to show they equal by realizing them as cycles arising in a graph
construction (\cite{MacPherson}, \cite{BFM}, or \cite{Fulton}, \S18.1).
The natural place to look for something of the sort is
$$\P(\bdll\oplus\bdlr)$$
There are two natural embeddings
$$\P(\bdll) \hookrightarrow \P(\bdll\oplus\bdlr)\quad, \quad \P(\bdlr)
\hookrightarrow \P(\bdll\oplus\bdlr)$$
as `first', resp.~`second' factor. These are the centers of two
{\it families of central projections}
$$\P(\bdll\oplus\bdlr) \overset\rho_{\Til M}\to\dashrightarrow
\P(\bdlr)\quad, \quad \P(\bdll\oplus\bdlr) \overset\rho_M\to
\dashrightarrow \P(\bdll)$$
respectively. Also, the rational maps from $\Til M$ to the bundles
(considered above) determine two embeddings
$$\gather
B\ell_{\pi^{-1}Y}\Til M \hookrightarrow \P(\bdll) \hookrightarrow
\P(\bdll\oplus\bdlr)\\
B\ell_{Y'}\Til M\hookrightarrow \P(\bdlr) \hookrightarrow
\P(\bdll\oplus\bdlr)\endgather$$
and the `cones'
$$G_M=\rho_M^{-1}(B\ell_{\pi^{-1}Y}\Til M)\quad,\quad G_{\Til
M}=\rho_{\Til M}^{-1}(B\ell_{Y'}\Til M)$$
(abusing notations). Also denote by $\Cal O(-1)$ the tautological
subbundle of $\P(\bdll\oplus\bdlr)$, and note that $\Cal O(-1)$
restricts to $\Cal O_M(-1)$, $\Cal O_{\Til M}(-1)$ on the two factors,
and that
$$\gather N_{\P(\bdll)}\P(\bdll\oplus\bdlr)=\bdlr\otimes\Cal O(1)\\
N_{\P(\bdlr)}\P(\bdll\oplus\bdlr)=\bdll\otimes\Cal O(1)
\endgather$$
(as seen with the aid of standard Euler sequences). Finally, $p$ will
denote the bundle map to $\Til M$. Here are some of the notations in a
diagram, for ease of reference:
$$\diagram
& \P(\bdll\oplus\bdlr) \xdashed[1,-1]_{\rho_M}|>\tip
\xdashed[1,1]^{\rho_{\Til M}}|>\tip \ddto^p \\
\P(\bdll) \drto_{p_M} & & \P(\bdlr) \dlto^{p_{\Til M}}\\
& \Til M \dto^\pi\\
& M
\enddiagram$$

Now we are ready for the new reformulation of what we have to prove:
\proclaim{Claim III.2} In order to prove (3), it suffices to show that
$$c\left(\frac{\bdll\oplus\bdlr}{\Cal O(-1)}\right) \cap
\left([G_{\Til M}]-[G_M]\right)$$
pushes forward to 0 in $M$.
\endproclaim
\demo{Proof} This follows immediately from Claim~III.1 and the following
Lemma:
\proclaim{Lemma III.1}
$$\align
p_*\left(c\left(\frac{\bdll\oplus\bdlr}{\Cal O(-1)}\right)
\cap[G_M]\right)&={p_M}_*\left(c\left(\frac{\bdll}{\Cal
O_M(-1)}\right)\cap [B\ell_{\pi^{-1}Y}\Til M]\right)\\
p_*\left(c\left(\frac{\bdll\oplus\bdlr}{\Cal O(-1)}
\right)\cap[G_{\Til M}]\right)&={p_{\Til M}}_*\left(c\left(
\frac{\bdlr}{\Cal O_{\Til M}(-1)}\right)\cap[B\ell_{Y'}\Til M]\right)
\endalign$$
\endproclaim
\demo{Proof} We check the first equality; the second is entirely
similar.

First, observe that since $G_M\cap \P(\bdll)=B\ell_{\pi^{-1}Y}\Til M$
in $\P(\bdll\oplus\bdlr)$ (and this intersection is transversal):
$${p_M}_*\left(c\left(\frac{\bdll}{\Cal
O_M(-1)}\right)\cap [B\ell_{\pi^{-1}Y}\Til M]\right) =
p_*\left(c\left(\frac{\bdll}{\Cal O_M(-1)}\right)\cap [\P(\bdll)]\cdot
[G_M]\right)$$
Next, $\Cal O(1)$ restricts to $\Cal O_M(1)$ on $\P(\bdll)$; denote by
$j$ its first Chern class. The normal bundle formula above tells
us that (with $n=\dim M$)
$$\align
[\P(\bdll)]\cdot [G_M] &=c_{top}(\bdlr\otimes \Cal O(1))\cap [G_M]\\
& =\left(j^{n+1}+\dots+j\,c_n(\bdlr)+c_{n+1}(\bdlr)\right)\cap [G_M]
\endalign$$
Therefore
$$\multline
c\left(\frac{\bdll}{\Cal O_M(-1)}\right)\cap [\P(\bdll)]\cdot [G_M]\\
=c(\bdll)(1+j+j^2+\dots)\left(j^{n+1}+\dots+j\,c_n(\bdlr)+c_{n+1}
(\bdlr)\right)\cap [G_M]
\endmultline$$
With the same notations:
$$\multline
c\left(\frac{\bdll\oplus\bdlr}{\Cal O(-1)}\right) \cap[G_M]\\
=c(\bdll)(1+j+j^2+\dots)\left(1+\dots+c_n(\bdlr)+c_{n+1}(\bdlr)\right)
\cap [G_M]
\endmultline$$
The difference consists of a sum of terms
$$c(\bdll)\,j^i c_k(\bdlr)\cap [G_M]$$
with $i+k<n+1$; but $G_M$ fibers over its image in $\P(\bdll)$ (via
$\rho_M$) with fibers of dimension $n+1$, so all such terms die
already in $\P(\bdll)$; and {\it a fortiori\/} after push forward to
$\Til M$.\qed\enddemo
\enddemo

\subhead \S3.3. (3) in terms of classes in $\P(\bdll\oplus\bdlr)$ over
$BL$\endsubhead
In order to attack Claim III.2, we have to produce explicitly a
class equivalent to $[G_{\Til M}]-[G_M]$, and we have to evaluate the
intersection product and push@-forward stated in Claim~III.2.
For this, we will pull@-back the situation to a variety dominating both
$B\ell_{\pi^{-1}Y} \Til M$ and $B\ell_{Y'} \Til M$. Consider the natural
morphism of bundles
$$\phi :\quad \bdll @>>> \bdlr$$
over $\Til M$, extending the differential $d\pi: T^*M @>>> T^*\Til M$.
(We are omitting here, and we will omit from now on, the pull@-back
notation $\pi^*$ on the sources of these morphisms.) This morphism will
play a fundamental role in what follows.

For a start, observe that $\phi$ is also a family of central
projections: over a general point of $\Til M$, $\phi$ is an
isomorphisms; over a point of $E$, say corresponding to a direction $u$
normal to $Z$, $\phi$ collapses forms vanishing along $TZ$ and $u$.
Projectivizing, we get a rational map
$$\psi:\quad \P\bdll \dashrightarrow \P\bdlr$$
which is resolved by blowing up the family $C$ of centers of the
projections (with the reduced structure); equivalently, the blow@-up
will be the graph $\Gamma$ of $\psi$ in $\P\bdll \times_{\Til M}
\P\bdlr$. At the same time, $\psi$ restricts to a rational map
$$B\ell_{\pi^{-1}Y}\Til M \dashrightarrow B\ell_{Y'}\Til M$$
which can be resolved by blowing up the source along its intersection
with $C$, obtaining a variety $BL$. Equivalently, $BL$ is the graph of
this map, which sits in $\Gamma$.
$$\diagram
& & \Gamma \ddllto
           \ddrrto\\
& & BL \uto \xto'[1,-1][2,-2] \xto'[1,1][2,2]\\
\P\bdll \xdashed[0,4]^{\psi}|>\tip & & & & \P\bdlr \\
B\ell_{\pi^{-1}Y}\Til M \uto \xdashed[0,4]|>\tip & & 
& & B\ell_{Y'}\Til M \uto
\enddiagram$$

Since $BL$ maps to both $B\ell_{\pi^{-1}Y}\Til M$, $B\ell_{Y'}\Til M$,
note that on $BL$ we have line bundles (obtained by pulling back) $\Cal
O_M(-1)=\Cal O(\Cal Y_M)$, $\Cal O_{\Til M}(-1)=\Cal O(\Cal Y_{\Til
M})$.

Now we can pull@-back $\P(\bdll\oplus\bdlr)$ etc.~to $BL$. The
advantage of doing so is that the cycles playing the role of $G_M$,
$G_{\Til M}$ have a nicer description: consider the following loci
defined over $BL$:
$$\gather
G_M=\P(\Cal O_M(-1)\oplus \bdlr)\subset \P(\bdll\oplus\bdlr)\quad\\
G_{\Til M}=\P(\bdll\oplus \Cal O_{\Til M}(-1)) \subset
\P(\bdll\oplus\bdlr)\quad;
\endgather$$
then the reader will verify that these $G_M$, $G_{\Til M}$ push forward
to the loci with the same name over $\Til M$ (the reason is that
$\P(\Cal O_M(-1))\subset \P(\bdll)$ realizes the embedding of
$B\ell_{\pi^{-1}Y}M$ in $\P\bdll$ etc.).

It follows that {\it we can adopt Claim III.2, taking the bundle and cycles
in the statement to live now over $BL$,\/} with the above positions.

\subhead \S3.4. The graph construction\endsubhead
The loci $G_M$, $G_{\Til M}$ are projectivizations of rank@-$(n+2)$
subbundles of $\bdll\oplus\bdlr$. It is natural to interpolate them by
considering the span of the graph of
$$\frac 1\lambda\phi:\quad\bdll @>>> \bdlr$$
and $0\oplus \Cal O_{\Til M}(-1)$ in $\bdll \oplus \bdlr$ : let
$G_\lambda$ denote this span (thus $G_\lambda$ is a rank@-$(n+2)$
subbundle of $\bdll\oplus\bdlr$ for all $\lambda\ne 0$).

Notice that as $1/\lambda \to 0$, $G_\lambda$ projectivizes to the
locus $G_{\Til M}$ defined above. In fact $G_{\lambda}$ can be realized
(for $\lambda\ne 0$) as the element in
$\Grass_{n+2}(\bdll\oplus\bdlr)$ determined by the graph of
$$\bdll @>{\frac 1\lambda \phi}>> \bdlr @>>> \frac \bdlr{\Cal
O_{\Til M}(-1)}\quad,\tag *$$
a point of $\Grass_{n+1}(\bdll\oplus (\bdlr/{\Cal O_{\Til M}(-1)}))$.
Over general points of $BL$, the ker of the composition (*) is the
fiber of $\Cal O_M(-1)$; it follows from general considerations
about the graph construction (\cite{BFM}; also \cite{Kwieci\'nski}, I.7,
p.~56) that the flat limit of $G_\lambda$ as $1/\lambda \to \infty$
will consist of several components, one of which will {\it precisely
projectivize to the locus\/} $G_M$ defined above. This will also be
recovered later on, by a coordinate computation.

Let $Z_\infty$ denote the {\it other\/} (that is, distinct from $G_M$)
components of the projectivization of $\lim_{1/\lambda \mapsto \infty}
G_{\lambda}$.
\proclaim{Claim III.3} In order to prove (3), it suffices to show that
$$c\left(\frac{\bdll\oplus\bdlr}{\Cal O(-1)}\right) \cap [Z_\infty]$$
pushes forward to 0 in $M$.
\endproclaim
\demo{Proof} The construction gives an explicit rational equivalence
between $[G_{\Til M}]=[\lim_{1/\lambda \mapsto 0} \P(G_\lambda)]$ and
$[\lim_{1/\lambda\mapsto\infty}\P(G_\lambda)]=[G_M]+[Z_\infty]$.
\qed\enddemo

\subhead \S3.5. Coordinate set@-up in $BL$\endsubhead
Fortunately $\lim_{\lambda \mapsto 0}G_\lambda$ (and therefore
$Z_\infty$) can be analyzed most explicitly by a coordinate computation;
however, this requires studying $BL$ more carefully and introducing if
possible yet more notations.

The variety $BL$ contains (inverse images of) divisors $E$, $\Cal Y_M$,
$\Cal Y_{\Til M}$; in fact we already observed that the line bundles
for the latter two are respectively $\Cal O_M(-1)$, $\Cal O_{\Til
M}(-1)$. Also, as $BL$ arises by blowing up $B\ell_{\pi^{-1}Y}\Til
M$ along a certain locus (called $C$ in the above), it contains a
corresponding exceptional divisor $\Cal E_M$. Further, note that the
rational map $\psi$ defined above is birational; so we can think of
$BL$ as the resolution of indeterminacies of $\psi^{-1}$, which
realizes it as a blow@-up of $B\ell_{Y'}\Til M$. Call $\Cal E_{\Til M}$
the exceptional divisor of this blow@-up. (Note: $\Cal E_M$, $\Cal
E_{\Til M}$ are restrictions of analogous exceptional divisors from
$\Gamma$.)

The intersection of $C$ and $B\ell_{\pi^{-1}Y}\Til M$ can be computed
easily. The reader may wish to check that it is the residual to $\Cal
Y_M$ of the scheme@-theoretic inverse image of $Y'$ in
$B\ell_{\pi^{-1}Y}\Til M$. Similarly, it is not hard to see that the
center of the blow@-up of $B\ell_{Y'}\Til M$ producing $BL$ is the
residual to $\Cal Y_{\Til M}$ in the (scheme theoretic) union of $E$
and the inverse image of $Y$ in $B\ell_{Y'}\Til M$.

To get a feeling for the situation, the reader may also wish to check that
the residual to $Y$ in $Y'$ is supported on the intersection of $E$ and the
proper transform $\Til X$ of $X$ in $\Til M$. In fact, $E\cap Y'$ is
precisely this intersection, while $E\cap Y$ is supported `just' on the
points at which $\Til X$ is tangent to the fibers of $E$ over $Z$.

\proclaim{Lemma III.2} With the above notations:
$$\gather
\Cal E_M+\Cal E_{\Til M}=E\\
\Cal Y_M+\Cal E_M=\Cal Y_{\Til M}
\endgather$$
as divisors of $BL$.\endproclaim
\demo{Proof} These equalities follow easily from the considerations
immediately preceding this statement. For example, the ideal of $\Cal
E_M$ is the pull@-back of the ideal of $C$, hence it is the residual to
$\Cal Y_M$ in $\Cal Y_{\Til M}$: this gives the second equality. By the
same token, $\Cal E_{\Til M}$ is the residual to $\Cal Y_{\Til M}$ in
$\Cal Y_M+E$: that is,
$$\Cal Y_M+E=\Cal Y_{\Til M}+\Cal E_{\Til M}$$
and the first equality follows.\qed\enddemo

Lemma III.2 will be used in a moment, when we choose functions on $BL$ 
to write entries for a matrix whose row@-span is $G_{\lambda}$. First,
we have to choose local coordinates in $\bdll$ and $\bdlr$. Choose local
parameters $x_1,\dots,x_n$ in $M$ and $\tilde x_1,\dots,\tilde x_n$ in
$\Til M$ so that $Z$ has ideal $(x_1,\dots,x_d)$ and the blow@-up map
$\Til M @>>> M$ is given by
$$\left\{\aligned x_1 &=\tilde x_1\\
x_2 &=\tilde x_1\tilde x_2\\
&\dots\\
x_d &=\tilde x_1\tilde x_d\\
x_{d+1}&=\tilde x_{d+1}\\
&\dots\\
x_n &=\tilde x_n
\endaligned\right.$$
Keeping in mind the sequences
$$\gather
0 @>>> \Cal T^*M\otimes \Cal L @>>> \bdll @>>> \Cal L @>>> 0\\
0 @>>> \Cal T^*\Til M\otimes \Cal L @>>> \bdlr @>>> \Cal L @>>> 0
\endgather$$
and working locally, we use $(s,v_1,\dots,v_n)$ to denote the jet in
$\bdll$ mapping to $s$ in the fibers of $\Cal L$ and with differential
$v_1 dx_1+\dots+v_n dx_n$. Similarly $(s,\tilde v_1,\dots,\tilde v_n)$
in $\bdlr$ maps to $s$ in $\Cal L$ and has differential $\tilde v_1
d\tilde x_1+\dots+\tilde v_n d\tilde x_n$. The morphism
$$\phi:\quad \bdll @>>> \bdlr$$
defined above has then matrix
$$\pmatrix
1 & 0 & 0 & 0 & \cdots & 0 & 0 & \cdots & 0 \\
0 & 1 &\tilde x_2 &\tilde x_3 & \cdots &\tilde x_d & 0 & \cdots & 0 \\ 
0 & 0 & \tilde x_1 & 0 & \cdots & 0 & 0 & \cdots & 0 \\
0 & 0 & 0 & \tilde x_1 & \cdots & 0 & 0 & \cdots & 0 \\
  & \vdots & & & \ddots & & & & \vdots\\
0 & 0 & 0 & 0 & \cdots & \tilde x_1 & 0 & \cdots & 0 \\
0 & 0 & 0 & 0 & \cdots & 0 & 1 & \cdots & 0 \\
  & \vdots & & & \vdots & & & \ddots & \vdots\\
0 & 0 & 0 & 0 & \cdots & 0 & 0 & \cdots & 1
\endpmatrix$$
in these coordinates. Also, the embeddings 
$$\gather
B\ell_{\pi^{-1}Y}\Til M\subset \bdll\\
B\ell_{Y'}\Til M\subset  \bdlr
\endgather$$ 
are obtained by projectivizing (and closing) the image of the sections
$$\gather
p \mapsto \left(F(p),{\ddfi 1}|_p,\dots,{\ddfi n}|_p\right)\\
p \mapsto \left(F(p),{\ddft 1}|_p,\dots,{\ddft n}|_p\right)\\
\endgather$$ 
Here $F$ is the section of $\Cal L$ giving the original hypersurface
$X$, pulled@-back to $\Til M$ to give the hypersurface $X'$. If $X$ has
multiplicity $m$ along $Z$, $F$ will be a multiple of $\tilde x_1^m$ on
$\Til M$.

With these notations and over points at which $F$ or some of its
partials $\ddft i$ do not vanish (that is, away from $Y'$), the
subbundle $G_\lambda\subset \bdll\oplus \bdlr$ defined above is
spanned by the $(n+1)$ rows of
$$\left(\matrix
\lambda & 0 & 0 & 0 & \cdots & 0 & 0 & \cdots & 0 \\
0 & \lambda & 0 & 0 & \cdots & 0 & 0 & \cdots & 0 \\
0 & 0 & \lambda & 0 & \cdots & 0 & 0 & \cdots & 0 \\
0 & 0 & 0 & \lambda & \cdots & 0 & 0 & \cdots & 0 \\
  & \vdots & & & \ddots & & & & \vdots\\
0 & 0 & 0 & 0 & \cdots & \lambda & 0 & \cdots & 0 \\
0 & 0 & 0 & 0 & \cdots & 0 & \lambda & \cdots & 0 \\
  & \vdots & & & \vdots & & & \ddots & \vdots\\
0 & 0 & 0 & 0 & \cdots & 0 & 0 & \cdots & \lambda
\endmatrix\right.\left|\matrix
1 & 0 & 0 & 0 & \cdots & 0 & 0 & \cdots & 0 \\
0 & 1 & 0 & 0 & \cdots & 0 & 0 & \cdots & 0 \\
0 & \tilde x_2 & \tilde x_1 & 0 & \cdots & 0 & 0 & \cdots & 0 \\
0 & \tilde x_3 & 0 & \tilde x_1 & \cdots & 0 & 0 & \cdots & 0 \\
  & \vdots & & & \ddots & & & & \vdots\\
0 & \tilde x_d & 0 & 0 & \cdots & \tilde x_1 & 0 & \cdots & 0 \\
0 & 0 & 0 & 0 & \cdots & 0 & 1 & \cdots & 0 \\
  & \vdots & & & \vdots & & & \ddots & \vdots\\
0 & 0 & 0 & 0 & \cdots & 0 & 0 & \cdots & 1
\endmatrix\right)$$
together with the row vector
$$\left(\matrix
0 & 0 & 0 & 0 & \cdots & 0 & 0 & \cdots & 0 
\endmatrix\right.\left|\matrix
F & \ddft 1 & \cdots & \ddft n\endmatrix\right)$$
(accounting for the $0\oplus \Cal O_{\Til M}(-1)$ factor). The reason
why we introduced the variety $BL$ above is to be able to extend this
description to points of $Y'$. Still working locally, give names
to the sections corresponding to the various divisors on $BL$. We have:

---$E$, with local generator $\tilde x_1$ (borrowing its name from
$\Til M$, from which the generator is pulled back);

---the exceptional divisors $\Cal Y_M$, $\Cal Y_{\Til M}$; the ideals of
these are the pull@-backs of $\Cal I_{\pi^{-1}Y}$, $\Cal I_{Y'}$, and we
will call local generators for these (principal) ideals $y_M$, $y_{\Til
M}$ respectively;

---the exceptional divisors $\Cal E_M$, $\Cal E_{\Til M}$; local
generators for these will be called $e_M$, $e_{\Til M}$ respectively.

Lemma III.2 gives the following relations among these terms:
$$\tilde x_1=e_M e_{\Til M}\quad,\quad y_{\Til M}=y_M e_M$$
(so that also $y_{\Til M}e_{\Til M}=y_M\tilde x_1$).

Since $\Cal I_{Y'}$ pulls back to $(y_{\Til M})$, there must be (local)
$a_0,\dots,a_n$ over $BL$ so that
$$F=a_0 y_{\Til M}\quad,\quad \ddft 1=a_1 y_{\Til M}\quad,\quad \dots
\quad,\quad\ddft n=a_n y_{\Til M}$$
By the same token, there must be $b_0,\dots,b_n$ such that
$$F=b_0 y_M\quad,\quad \ddi F1=b_1 y_M\quad,\quad \dots \quad, \quad
\ddi Fn=b_n y_M$$
Further, both ideals $(a_0,\dots,a_n)$ and $(b_0,\dots,b_n)$ equal
$(1)$.

Now
$$\gather
\ddft 1=\ddi F1+\tilde x_2 \ddi F2+\dots+\tilde x_d \ddi Fd\\
\ddft 2=\tilde x_1 \ddi F2,\quad\dots\quad, \ddft d=\tilde
x_1 \ddi Fd\\
\ddft {d+1}=\ddi F{d+1},\quad\dots\quad, \ddft n=\ddi Fn
\endgather$$
which translates into
$$\gather
a_0 y_{\Til M}=b_0 y_M,\quad a_1 y_{\Til M}=b_1 y_M+\tilde x_2 b_2
y_M+\dots+\tilde y_d b_d y_M\\
a_2 y_{\Til M}=\tilde x_1 b_2 y_M,\quad\dots\quad, a_d y_{\Til M}=\tilde
x_1 b_d y_M\\
a_{d+1} y_{\Til M}=b_{d+1} y_M,\quad\dots\quad,a_n y_{\Til M}= b_n y_M
\endgather$$
and therefore
$$\gather
b_0=a_0 e_M,\quad b_1=a_1 e_M-\tilde x_2 b_2-\dots-\tilde x_d
b_d\\ a_2=b_2 e_{\Til M},\quad\dots\quad,a_d=b_d e_{\Til M}\\
b_{d+1}=a_{d+1} e_M,\quad\dots\quad, b_n=a_n e_M
\endgather$$
This means that we can throw away half of the $a$'s and $b$'s; from the
data of
$$\{a_0, a_1, b_2,\dots, b_d,a_{d+1},\dots,a_n\}$$
(which again locally generate $(1)$) we can obtain explicit coordinates
in $\P\bdll$, $\P\bdlr$ for the image of a point of $BL$:
$$\gather
(a_0 e_M:a_1 e_M-\tilde x_2 b_2-\dots-\tilde x_d
b_d:b_2:\dots: b_d:a_{d+1} e_M:\dots:a_n e_M),\\
(a_0:a_1:b_2e_{\Til M}:\dots: b_d e_{\Til M}:a_{d+1}:\dots:a_n)
\endgather$$

With this understood, the rows of the $(n+2)\times (2n+2)$@-matrix
$$\left(\matrix
\lambda & 0 & 0 & \cdots & 0 & 0 & \cdots & 0 \\
0 & \lambda & 0 & \cdots & 0 & 0 & \cdots & 0 \\
0 & 0 & \lambda & \cdots & 0 & 0 & \cdots & 0 \\
  & \vdots & & \ddots & & & & \vdots\\
0 & 0 & 0 & \cdots & \lambda & 0 & \cdots & 0 \\
0 & 0 & 0 & \cdots & 0 & \lambda & \cdots & 0 \\
  & \vdots & & \vdots & & & \ddots & \vdots\\
0 & 0 & 0 & \cdots & 0 & 0 & \cdots & \lambda\\
0 & 0 & 0 & \cdots & 0 & 0 & \cdots & 0 
\endmatrix\right.\left|\matrix
1 & 0 & 0  & \cdots & 0 & 0 & \cdots & 0 \\
0 & 1 & 0  & \cdots & 0 & 0 & \cdots & 0 \\
0 & \tilde x_2 & \tilde x_1 & \cdots & 0 & 0 & \cdots & 0 \\
  & \vdots & & \ddots & & & & \vdots\\
0 & \tilde x_d & 0 & \cdots & \tilde x_1 & 0 & \cdots & 0 \\
0 & 0 & 0 & \cdots & 0 & 1 & \cdots & 0 \\
  & \vdots & & \vdots & & & \ddots & \vdots\\
0 & 0 & 0 & \cdots & 0 & 0 & \cdots & 1\\
a_0 & a_1 & b_2 e_{\Til M} & \cdots & b_d e_{\Til M} & a_{d+1} &
\cdots & a_n
\endmatrix\right)$$
do span $G_\lambda$ for all $\lambda\ne 0$.

\subhead \S3.6. Three lemmas\endsubhead
The following lemmas will not be used till the final step of our proof;
but this seems the best place to include them, as they use the
notations we just introduced. Consider again the coordinates for a
point in $B\ell_{\pi^{-1}Y}\Til M\subset\P\bdll$ obtained above:
$$\left(\boxed{a_0 e_M}:a_1 e_M-\tilde x_2 b_2-\dots-\tilde x_d
b_d:b_2:\dots: b_d:\boxed{a_{d+1} e_M:\dots:a_n e_M}\right)$$
The entries here are (local) components of a vector spanning $\Cal
O_M(-1)$ in $\bdll$. It is clear that the boxed entries vanish along
$\Cal E_M$ (as $e_M=0$ is an equation for the latter). To express this
more intrinsically, consider the natural projection
$$\rho_Z:\quad\bdll\oplus \bdlr @>>> \bdll @>>> \Cal P^1_Z\Cal L$$
killing the second factor and projecting the first onto the bundle of
principal parts of $\Cal L$ over $Z$; then the above observation is that
\proclaim{Lemma III.3} $\rho_Z(\Cal O_M(-1)\oplus 0)=0$ along $\Cal E_M$.
\endproclaim

This vanishing will be an important ingredient at the final step of the
proof, in \S 3.8. The two lemmas that follow are also important, and they
are less evident. First, observe that as $\Cal E_{\Til M}$ maps
injectively into $B\ell_{\pi^{-1}}\Til M$, its components either 

(i) dominate $E\subset \Til M$, or 

(ii) dominate components in $E\subset B\ell_{\pi^{-1}}\Til M$ which
contract in $\Til M$.

\proclaim{Lemma III.4} $\rho_Z(\Cal O_M(-1)\oplus 0)=0$ along components of
$\Cal E_{\Til M}$ of type (i).\endproclaim

\demo{Proof} To see this, it suffices to check that, at points of $\Cal
E_{\Til M}$ mapping to general points of $E\subset\Til M$, necessarily
$e_{\Til M}$ divides $a_0$ and $a_{d+1},\dots,a_n$. Assume $X$ has
multiplicity $m\ge 1$ along $Z$, and write (locally) $F=\tilde x_1^m\Til
F$ (that is, let $(\Til F)$ be the ideal of the proper transform $\Til X$
of $X$). Then computing partials gives the ideal of $\pi^{-1}Y$:
$$\Cal I_{\pi^{-1}Y}=\tilde x_1^{m-1}\left(\tilde x_1 \Til F, m\Til
F+\tilde x_1\ddti F1,\ddti F2,\dots,\ddti Fd,\tilde x_1 \ddti F{d+1},
\dots,\tilde x_1\ddti Fn\right)\tag *$$
This ideal pulls back to $(y_M)$ in $BL$. Now $\Til F\ne 0$ at a
general point of $E\subset \Til M$, and hence at a general point $p$
of the component of $\Cal E_{\Til M}$ we are considering. Then at such
points (*) gives 
$$\Cal I_{\pi^{-1}Y}=(\tilde x_1^{m-1})\quad:$$
that is, $y_M=\tilde x_1^{m-1}$ in $BL$; with the positions made in
\S3.5:
$$a_0 e_M=\tilde x_1\Til F, a_{d+1} e_M=\tilde x_1 \ddti F{d+1}, \dots,
a_n e_M=\tilde x_1\ddti Fn\quad,$$
that is
$$a_0 =\Til F e_{\Til M}, a_{d+1}=\ddti F{d+1} e_{\Til M}, \dots,
a_n=\ddti Fn e_{\Til M}\quad:$$
$a_0$ and $a_{d+1},\dots,a_n$ are multiples of $e_{\Til M}$ near $p$, as
we needed.\qed\enddemo

The story for components of type (ii) is a little different: the boxed
entries above do {\it not} vanish identically along these components.
However, let $\Cal E'_{\Til M}$ be such a component, and let $Z'$ be
the subvariety of $Z$ that $\Cal E'_{\Til M}$ dominates; also, let
$$\bdll\oplus \bdlr @>\rho_{Z'}>> \Cal P^1_{Z'}\Cal L$$
denote the natural projection. Then:
\proclaim{Lemma III.5} $\rho_{Z'}(\Cal O_M(-1)\oplus 0)=0$ along $\Cal
E'_{\Til M}$.\endproclaim
\demo{Proof} Again, it suffices to show this at the general point $p$ of
$\Cal E'_{\Til M}$. Let then $\tilde y_M$ denote a local equation for
$\Cal E'_{\Til M}$ at $p$. $\Cal E'_{\Til M}$ dominates a component
contained in $E\subset B\ell_{\pi^{-1}}\Til M$ of the exceptional
divisor of the blow@-up of $\Til M$ along $\pi^{-1} Y$. Thus $\Cal
E'_{\Til M}$ is, aside of a multiple of $E$, the inverse image of
$\pi^{-1} Y$ in $BL$ (near $p$), and we may assume that the ideal of
the latter (again, aside of the factor $\tilde x_1^{m-1}$) contains
$\tilde x_1$: from (*) above, we see that $(\tilde y_M)$ is the
pull@-back of
$$\Cal J=\left(\tilde x_1,\Til F,\ddti F2,\dots,\ddti Fd\right)$$
near $p$, and (since $\tilde x_1^{m-1}\Cal J=\Cal I_{\pi^{-1}Y}$ pulls
back to $(y_M)$) we have
$$y_M=\tilde x_1^{m-1} \tilde y_M$$
up to units at $p$.

Since $\tilde x_1\in\Cal J$ and $\Til F\in \Cal J$, we have
$$\tilde x_1=c_x\tilde y_M \quad,\quad \Til F=c_F\tilde y_M$$
for some $c_x, c_F$. By the positions made in \S 3.5:
$$y_M a_0 e_M=F=\tilde x_1^m \Til F\quad,$$
that is
$$\tilde y_M a_0 e_M=\tilde x_1 \Til F\quad,$$
and therefore
$$a_0 e_M = c_x c_F \tilde y_M\quad.$$
The left@-hand@-side is the first boxed entry listed at the beginning
of this subsection, and the right@-hand@-side shows that this vanishes
along $\Cal E'_{\Til M}$ (as $\tilde y_M=0$ is an equation for the
latter). Behind the notational smoke, the reader should be able to see
that this simply works because both $\tilde x_1$ and $\Til F$ vanish
along the subscheme defined by $\Cal J$ in $\Til M$; the first
boxed entry is controlled by $\tilde x_1\Til F$, so it vanishes to
higher order. The above computation simply formalizes this observation.

Now we want to argue similarly for the other entries. For
$i=d+1,\dots,n$ we have (from \S 3.5)
$$y_M a_i e_M= \ddft i =\tilde x_1^m\ddti Fi$$
and therefore
$$\tilde y_M a_i e_M=\tilde x_1\ddti Fi$$
or
$$a_i e_M= c_x\ddti Fi\quad;$$
the vanishing of $a_i e_M$ along $\Cal E'_{\Til M}$ follows for $i>d+1$
if $\ddti Fi$ vanishes along the subset $S\subset \Til M$ dominated by
$\Cal E'_{\Til M}$.

Recall we are denoting by $Z'\subset Z$ the image of $\Cal E'_{\Til
M}$ (hence of $S$). The image of $p$ will be a general, hence
nonsingular point $q$ of $Z'$. Choose the local parameters $\tilde
x_{d+1},\dots,\tilde x_n$ on $Z$ so that $Z'$ has equations $\tilde
x_{d+1}=\dots=\tilde x_s=0$ in $Z$ near $q$. Now $Z'$ is the image of
$S\subset \Til F$: the (Zariski) tangent space to $\Til F$ at a
general point of $S$ must dominate the tangent space to $Z'$. This
means that $\ddti Fi=0$ along $S$ for $i=s+1,\dots,n$, and this implies
the vanishing of the corresponding entries $a_i e_M$ along $\Cal E'_{\Til
M}$, as observed above. Intrinsically, this amounts to the vanishing of
$\rho_{Z'}(\Cal O_M(-1)\oplus 0)$, and we are done.\qed\enddemo

\example{Example} The case considered in Lemma~III.5 occurs when the
hypersurface is nonsingular away from $Z$, and `cuspidal' along some
subset $Z'$ of $Z$. For an example that may help fixing ideas,
consider the surface with equation
$$x_2^2+x_1^2(x_1+x_3)=0$$
in $\Bbb A^3$. This is singular (and equimultiple) along the line
$x_1=x_2=0$; the singularity has transversal branches at all points
with $x_3\ne 0$, but is cuspy at the origin. In the blow@-up, the
proper transform of the surface intersects each fiber of the
exceptional divisor in two points; these collide into one point,
say $r$, in the fiber over the origin. In the above terminology, $S=\{r\}$;
so blowing up $\pi^{-1}Y$ in $\Til M$ amounts to blowing up $r$ in this
case. Since $r\in E$, this will produce a component of $E\subset
B\ell_{\pi^{-1}Y}\Til M$ contracting into $\Til M$, that is a component of
type (ii) in $\Cal E_{\Til M}$. In this example $Z'$ would be the
origin.\endexample

\subhead \S3.7. Computing $Z_\infty$\endsubhead
In order to determine the limit of $G_\lambda$ as $\lambda\mapsto 0$ we
consider the matrix given in \S3.5 as defining a rational map
$BL\times \P^1 \dashrightarrow \Grass_{n+2}(\bdll\oplus \bdlr)$, by
sending (for $\lambda\ne 0$) $(p,(\lambda:1))$ to the fiber of
$G_{\lambda}$ over $p$. The plan is to resolve the indeterminacies of
this map, and determine $Z_\infty$ as the image of specific loci via the
resolved map.

The base locus of the map is determined (thinking of Pl\"ucker
coordinates for $\Grass_{n+2}$) by the ideal of $(n+2)\times (n+2)$
minors of the matrix of row@-vectors given in \S3.5. In terms of
the $\ddft i$'s, this turns out to be the ideal
$$\align
& \lambda^d \left(F,\ddft 1,\dots,\ddft n\right)\\
+&\lambda^{d-1} \left(\tilde x_1 F,\tilde x_1 \ddft1, \ddft2,\dots,
\ddft d, \tilde x_1 \ddft{d+1},\dots,\tilde x_1 \ddft n\right)\\
+&\dots\\
+&\lambda \tilde x^{d-2} \left(\tilde x_1 F,\tilde x_1 \ddft1,
\ddft2,\dots, \ddft d, \tilde x_1 \ddft{d+1},\dots,\tilde x_1 \ddft
n\right)
\endalign$$
Pulling back to $BL$ this is written
$$\lambda y_{\Til M}\left(\lambda^{d-1},\lambda^{d-2} e_{\Til M},
\lambda^{d-3}e_Me_{\Til M}^2,\dots,e_M^{d-2}e_{\Til M}^{d-1}\right)$$
and resolving the map amounts to making the ideal in $(\,)$
principal. 

By our good fortune, this is easy to accomplish: it suffices to blow@-up
$BL\times \P^1$ twice, first along $\Cal E_{\Til M}\subset \Lambda$,
where $\Lambda$ is the copy $BL\times \{(0:1)\}$ of $BL$, and then along
the proper transform of $\Cal E_M\subset \Lambda$. In terms
of ideals, first we blow@-up along $(\lambda,e_{\Til M})$; the
interesting chart (we leave to the reader checking that nothing goes
wrong on the other charts in our blow@-ups) is
$$\left\{\aligned\lambda &=\tilde \lambda \tilde
e_{\Til M}\\ e_{\Til M} &=\tilde e_{\Til M}\endaligned\right.\quad;$$
so the ideal pulls back to
$$\left(\tilde \lambda^{d-1} \tilde e_{\Til M}^{d-1}, \tilde
\lambda^{d-2} \tilde e_{\Til M}^{d-1}, \tilde \lambda^{d-3} \tilde
e_{\Til M}^{d-1} e_M,\dots, \tilde e_{\Til M}^{d-1}
e_M^{d-2}\right)=\tilde e_{\Til M}^{d-1}(\tilde \lambda,e_M)^{d-2}
\quad;$$
then along $(\tilde \lambda,e_M)$. It is clear the pull@-back of the
ideal will then be principal, as claimed.

Now we have to do this on the matrix whose rows span the
$G_\lambda$'s. Again, we show what happens on one interesting chart of
the result, and leave the others to the reader. We change coordinates
according to
$$\left\{\aligned \lambda&=s t\\
e_{\Til M} &=t\endaligned\right.\qquad,
\qquad \left\{\aligned s&=u v\\
e_M&=v\endaligned\right.$$
that is
$$\left\{\aligned \lambda&=uvt\\ e_M&=v \\ e_{\Til M} &=t \\ \tilde x_1
&=vt \endaligned\right.\quad,$$
and the matrix becomes{\eightpoint
$$\left(\matrix
uvt & 0 & 0 & \cdots & 0 & 0 & \cdots & 0 \\
0 & uvt & 0 & \cdots & 0 & 0 & \cdots & 0 \\
0 & 0 & uvt & \cdots & 0 & 0 & \cdots & 0 \\
  & \vdots & & \ddots & & & & \vdots\\
0 & 0 & 0 & \cdots & uvt & 0 & \cdots & 0 \\
0 & 0 & 0 & \cdots & 0 & uvt & \cdots & 0 \\
  & \vdots & & \vdots & & & \ddots & \vdots\\
0 & 0 & 0 & \cdots & 0 & 0 & \cdots & uvt\\
0 & 0 & 0 & \cdots & 0 & 0 & \cdots & 0
\endmatrix\right.\left|\matrix
1 & 0 & 0 & \cdots & 0 & 0 & \cdots & 0 \\
0 & 1 & 0 & \cdots & 0 & 0 & \cdots & 0 \\
0 & \tilde x_2 & vt & \cdots & 0 & 0 & \cdots & 0 \\
  & \vdots & & \ddots & & & & \vdots\\
0 & \tilde x_d & 0 & \cdots & vt & 0 & \cdots & 0 \\
0 & 0 & 0 & \cdots & 0 & 1 & \cdots & 0 \\
  & \vdots & & \vdots & & & \ddots & \vdots\\
0 & 0 & 0 & \cdots & 0 & 0 & \cdots & 1\\
a_0 & a_1 & b_2 t & \cdots & b_d t & a_{d+1} & \cdots & a_n
\endmatrix\right)$$}
At general points, this has the same span as{\eightpoint
$$\left(\matrix
uvt & 0 & 0 & \cdots & 0 & 0 & \cdots & 0 \\
0 & uvt & 0 & \cdots & 0 & 0 & \cdots & 0 \\
0 & -\tilde x_2 u & u & \cdots & 0 & 0 & \cdots & 0 \\
  & \vdots & & \ddots & & & & \vdots\\
0 & -\tilde x_d u & 0 & \cdots & u & 0 & \cdots & 0 \\
0 & 0 & 0 & \cdots & 0 & uvt & \cdots & 0 \\
  & \vdots & & \vdots & & & \ddots & \vdots\\
0 & 0 & 0 & \cdots & 0 & 0 & \cdots & uvt\\
0 & 0 & 0 & \cdots & 0 & 0 & \cdots & 0
\endmatrix\right.\left|\matrix
1 & 0 & 0 & \cdots & 0 & 0 & \cdots & 0 \\
0 & 1 & 0 & \cdots & 0 & 0 & \cdots & 0 \\
0 & 0 & 1 & \cdots & 0 & 0 & \cdots & 0 \\
  & \vdots & & \ddots & & & & \vdots\\
0 & 0 & 0 & \cdots & 1 & 0 & \cdots & 0 \\
0 & 0 & 0 & \cdots & 0 & 1 & \cdots & 0 \\
  & \vdots & & \vdots & & & \ddots & \vdots\\
0 & 0 & 0 & \cdots & 0 & 0 & \cdots & 1\\
a_0 & a_1 & b_2 t & \cdots & b_d t & a_{d+1} & \cdots & a_n
\endmatrix\right)$$}
\noindent hence as{\eightpoint
$$\left(\matrix
uvt & 0 & 0 & \cdots & 0 & 0 & \cdots & 0 \\
0 & uvt & 0 & \cdots & 0 & 0 & \cdots & 0 \\
0 & -\tilde x_2 u & u & \cdots & 0 & 0 & \cdots & 0 \\
  & \vdots & & \ddots & & & & \vdots\\
0 & -\tilde x_d u & 0 & \cdots & u & 0 & \cdots & 0 \\
0 & 0 & 0 & \cdots & 0 & uvt & \cdots & 0 \\
  & \vdots & & \vdots & & & \ddots & \vdots\\
0 & 0 & 0 & \cdots & 0 & 0 & \cdots & uvt\\
a_0 v & a_1 v-\sum_2^d \tilde x_i b_i & b_2 & \cdots & b_d & a_{d+1}
v & \cdots & a_n v\endmatrix\right.\left|\matrix
1 & 0 & 0 & \cdots & 0 & 0 & \cdots & 0 \\
0 & 1 & 0 & \cdots & 0 & 0 & \cdots & 0 \\
0 & 0 & 1 & \cdots & 0 & 0 & \cdots & 0 \\
  & \vdots & & \ddots & & & & \vdots\\
0 & 0 & 0 & \cdots & 1 & 0 & \cdots & 0 \\
0 & 0 & 0 & \cdots & 0 & 1 & \cdots & 0 \\
  & \vdots & & \vdots & & & \ddots & \vdots\\
0 & 0 & 0 & \cdots & 0 & 0 & \cdots & 1\\
0 & 0 & 0 & \cdots & 0 & 0 & \cdots & 0
\endmatrix\right)$$}
This matrix has maximal rank everywhere, as it ought to (this is
because $(a_0 v,a_1 v-\sum_2^d \tilde x_i b_i,b_2,\cdots,b_d,a_{d+1} v,
\cdots, a_n v)=(1)$, cf.~\S 3.5); so the corresponding map to
$\Grass_{n+2}$ is indeed defined everywhere. The $\lim_{\lambda\mapsto
0}G_\lambda$ is therefore determined by the image of $\lambda=0$; since
$\lambda=uvt$, this breaks up the limit into three pieces:

---over $u=0$, that is the component dominating $BL\times\{(0:1)\}$,
this gives{\eightpoint
$$\left(\matrix 0 & 0 & 0 & \cdots & 0 & 0 & \cdots & 0 \\
0 & 0 & 0 & \cdots & 0 & 0 & \cdots & 0 \\
 0 & 0 & 0 & \cdots & 0 & 0 & \cdots & 0 \\
  & \vdots & & \ddots & & & & \vdots\\
0 & 0 & 0 & \cdots & 0 & 0 & \cdots & 0 \\
0 & 0 & 0 & \cdots & 0 & 0 & \cdots & 0 \\
  & \vdots & & \vdots & & & \ddots & \vdots\\
0 & 0 & 0 & \cdots & 0 & 0 & \cdots & 0\\
a_0 v & a_1 v-\sum_2^d \tilde x_i b_i & b_2 & \cdots & b_d & a_{d+1} v
& \cdots & a_n v\endmatrix\right.\left|\matrix
1 & 0 & 0 & \cdots & 0 & 0 & \cdots & 0 \\
0 & 1 & 0 & \cdots & 0 & 0 & \cdots & 0 \\
0 & 0 & 1 & \cdots & 0 & 0 & \cdots & 0 \\
  & \vdots & & \ddots & & & & \vdots\\
0 & 0 & 0 & \cdots & 1 & 0 & \cdots & 0 \\
0 & 0 & 0 & \cdots & 0 & 1 & \cdots & 0 \\
  & \vdots & & \vdots & & & \ddots & \vdots\\
0 & 0 & 0 & \cdots & 0 & 0 & \cdots & 1\\
0 & 0 & 0 & \cdots & 0 & 0 & \cdots & 0
\endmatrix\right)$$}
which projectivizes to $[\P(\Cal O_M(-1)\oplus\bdlr)]=[G_M]$ (as we
promised in \S3.4, $G_M$ had to appear as one component in the limit);

---over $v=0$, which dominates $\Cal E_M\times \{(0:1)\}$:{\eightpoint
$$\left(\matrix
0 & 0 & 0 & \cdots & 0 & 0 & \cdots & 0 \\
0 & 0 & 0 & \cdots & 0 & 0 & \cdots & 0 \\
0 & -\tilde x_2 u & u & \cdots & 0 & 0 & \cdots & 0 \\
  & \vdots & & \ddots & & & & \vdots\\
0 & -\tilde x_d u & 0 & \cdots & u & 0 & \cdots & 0 \\
0 & 0 & 0 & \cdots & 0 & 0 & \cdots & 0 \\
  & \vdots & & \vdots & & & \ddots & \vdots\\
0 & 0 & 0 & \cdots & 0 & 0 & \cdots & 0\\
0 & -\sum_2^d \tilde x_i b_i & b_2 & \cdots & b_d & 0
& \cdots & 0\endmatrix\right.\left|\matrix
1 & 0 & 0 & \cdots & 0 & 0 & \cdots & 0 \\
0 & 1 & 0 & \cdots & 0 & 0 & \cdots & 0 \\
0 & 0 & 1 & \cdots & 0 & 0 & \cdots & 0 \\
  & \vdots & & \ddots & & & & \vdots\\
0 & 0 & 0 & \cdots & 1 & 0 & \cdots & 0 \\
0 & 0 & 0 & \cdots & 0 & 1 & \cdots & 0 \\
  & \vdots & & \vdots & & & \ddots & \vdots\\
0 & 0 & 0 & \cdots & 0 & 0 & \cdots & 1\\
0 & 0 & 0 & \cdots & 0 & 0 & \cdots & 0
\endmatrix\right)$$}

---and over $t=0$, which dominates $\Cal E_{\Til M}\times
\{(0:1)\}$:{\eightpoint
$$\left(\matrix
0 & 0 & 0 & \cdots & 0 & 0 & \cdots & 0 \\
0 & 0 & 0 & \cdots & 0 & 0 & \cdots & 0 \\
0 & -\tilde x_2 u & u & \cdots & 0 & 0 & \cdots & 0 \\
  & \vdots & & \ddots & & & & \vdots\\
0 & -\tilde x_d u & 0 & \cdots & u & 0 & \cdots & 0 \\
0 & 0 & 0 & \cdots & 0 & 0 & \cdots & 0 \\
  & \vdots & & \vdots & & & \ddots & \vdots\\
0 & 0 & 0 & \cdots & 0 & 0 & \cdots & 0\\
a_0 v & a_1 v-\sum_2^d \tilde x_i b_i & b_2 & \cdots & b_d & a_{d+1} v
& \cdots & a_n v\endmatrix\right.\left|\matrix
1 & 0 & 0 & \cdots & 0 & 0 & \cdots & 0 \\
0 & 1 & 0 & \cdots & 0 & 0 & \cdots & 0 \\
0 & 0 & 1 & \cdots & 0 & 0 & \cdots & 0 \\
  & \vdots & & \ddots & & & & \vdots\\
0 & 0 & 0 & \cdots & 1 & 0 & \cdots & 0 \\
0 & 0 & 0 & \cdots & 0 & 1 & \cdots & 0 \\
  & \vdots & & \vdots & & & \ddots & \vdots\\
0 & 0 & 0 & \cdots & 0 & 0 & \cdots & 1\\
0 & 0 & 0 & \cdots & 0 & 0 & \cdots & 0
\endmatrix\right)$$}
These last two loci make up $Z_\infty$ (by definition of the latter as
residual of $G_M$ in the limit). We have to study these loci, aiming
toward computing the class of the statement of Claim~III.3.

\subhead \S3.8. End of the proof of (3)\endsubhead
Summarizing, we have determined $Z_\infty$ as the image of two loci
defined over a double blow@-up of $BL\times\P^1$. They both sit in (the
 pull@-back of) $\P(\bdll\oplus\bdlr)$, and they dominate respectively
the second and first exceptional divisors (equations $v=0$, $t=0$),
which we will name $\Cal D_M$, $\Cal D_{\Til M}$; these in turn dominate
resp.~$\Cal E_M$, $\Cal E_{\Til M}$. The loci can be written as
projectivizations $F_M=\P\Cal F_M$, $F_{\Til M}=\P\Cal F_{\Til M}$ of
rank@-$(n+2)$ subbundles $\Cal F_M$, $\Cal F_{\Til M}$ of
$\bdll\oplus\bdlr$ over resp.~$\Cal D_M$, $\Cal D_{\Til M}$, determined by
the last two matrices written above. The following is our final
reformulation of (3):
\proclaim{Claim III.4} In order to prove (3), it suffices to show that
$$c\left(\frac{\bdll\oplus\bdlr}{\Cal O(-1)}\right)\cap [F_M]\quad,
\quad c\left(\frac{\bdll\oplus\bdlr}{\Cal O(-1)}\right)\cap [F_{\Til
M}]$$
vanish after push@-forward to $M$.\endproclaim
\demo{Proof} Via the map of bundles
$$\CD
\P(\bdll\oplus\bdlr) @>>> \P(\bdll\oplus\bdlr)\\
@VVV @VVV\\
\text{double blow@-up of $BL\times\P^1$} @>>> BL
\endCD$$
the cycle $[F_M]+[F_{\Til M}]$ (which lives in the top@-left spot)
pushes forward to $[Z_\infty]$ (in the top@-right spot). So the claim
follows directly from Claim~III.3.\qed\enddemo

Finally, we are ready to complete the proof of (3), and therefore of
the main Theorem:
\demo{Proof of (3)} Observe that
$$c\left(\frac{\bdll\oplus\bdlr}{\Cal O(-1)}\right)\cap [F_M] =
c\left(\frac{\bdll\oplus\bdlr}{\Cal F_M}\right)\, c\left(\frac{\Cal
F_M}{\Cal O(-1)}\right)\cap [F_M]$$
pushes forward to
$$c\left(\frac{\bdll\oplus\bdlr}{\Cal F_M}\right)\cap [\Cal D_M]
\tag{$\dagger$}$$
on the double blow@-up of $BL\times\P^1$: indeed, $\Cal O(-1)$ restricts
to the universal line bundle in $\Cal F_M$, so $c(\Cal O(-1))^{-1}\cap
[F_M]$ pushes forward to $c(\Cal F_M)^{-1}\cap [\Cal D_M]$. Similarly
$$c\left(\frac{\bdll\oplus\bdlr}{\Cal O(-1)}\right)\cap [F_{\Til M}]$$
pushes forward to
$$c\left(\frac{\bdll\oplus\bdlr}{\Cal F_{\Til M}}\right)\cap [\Cal
D_{\Til M}]\tag{$\dagger\dagger$}$$

The reason why these classes ($\dagger$), ($\dagger\dagger$) vanish
when pushed forward to $M$ lies in the three Lemmas in \S3.6. Arguing
explicitly for ($\dagger$), the key observation is that $\Cal F_M$ is
contained in the kernel of the natural morphism $\bdll \oplus \bdlr
@>\rho_Z>> \Cal P^1_Z \Cal L$: this can be checked locally, and it is
immediate from the matrix description given above, since in the chosen
coordinates $\rho_Z$ acts
$$\left(\matrix v_0 & \dots & v_n\endmatrix\right. \left|\matrix
\tilde v_0 & \cdots & \tilde v_n\endmatrix\right) \mapsto\left(\matrix
v_0 & v_{d+1} & \cdots & v_n\endmatrix\right)$$
Lemma~III.3 in \S3.6 amounts to observing this vanishing for the last
row of the matrix, as $\Cal D_M$ dominates $\Cal E_M$; the vanishing
for the rest of the matrix is clear for ($\dagger$) as well as for
($\dagger\dagger$).

Therefore we have an onto morphism
$$\frac{\bdll\oplus\bdlr}{\Cal F_M} @>>> \Cal P^1_Z\Cal L @>>>
0\quad;$$
if $\Cal K$ denotes the kernel of this morphism, and $\Pi$ denotes the
projection to $M$, we get
$$\Pi_*\left(c\left(\frac{\bdll\oplus\bdlr}{\Cal F_M}\right)\cap [\Cal
D_M]\right)=c(\Cal P^1_Z\Cal L) \Pi_*(c(\Cal K)\cap [\Cal D_M])\quad.$$
Now $\Cal D_M$ has dimension $n$ ($\Cal D_M$ is a divisor in a
blow@-up of $BL\times\P^1$, and $BL$ is birational to $M$) while $\Cal
K$ has rank $(2n+2)-(n+2)-(n-d+1)=d-1$, so that $c(\Cal K)\cap [\Cal
D_M]$ has no terms in dimension $\le\dim Z$; while $\Cal D_M$
dominates $Z$ via $\Pi$: $\Cal D_M$ dominates $\Cal E_M$, then $E$,
then $Z$. This forces the last $\Pi_*$ to vanish, as needed.

Concerning ($\dagger\dagger$), Lemma~III.4 in \S3.6 shows that $\Cal F_{\Til
M}$ is in the kernel of $\rho_Z$ along components dominating
components `of type (i)' of $\Cal E_{\Til M}$, and the vanishing follows
by the same argument as for ($\dagger$).

The situation is slightly more complicated over components `of type
(ii)'. By Lemma~III.5 in \S3.6 we know that, along such a component $\Cal
D'_{\Til M}$, $\Cal F_{\Til M}$ is in the kernel of the epimorphism
$\rho_{Z'}: \bdll\oplus\bdlr @>>> \Cal P^1_{Z'}\Cal L$ for the
subvariety $Z'$ of $Z$ dominated by the component. Now pull@-back the
situation through the fiber square
$$\CD
\Til{\Cal D}'_{\Til M} @>>> {\Cal D}'_{\Til M}\\
@VVV @VVV \\
\Til Z' @>>> Z'
\endCD$$
where the bottom row is the Nash@-blow@-up of $Z'$: over $\Til Z'$, the
pull@-back of $\Cal P^1_{Z'}\Cal L$ surjects onto a locally free sheaf
of rank $1+\dim Z'$; the above argument then gives that
$$\left(c\left(\frac{\bdll\oplus\bdlr}{\Cal F_{\Til M}}\right)\cap [
\Til{\Cal D}'_{\Til M}]\right)$$
vanishes after push@-forward to $\Til Z'$, and this implies the
vanishing of
$$c\left(\frac{\bdll\oplus\bdlr}{\Cal F_{\Til M}}\right)\cap [\Cal
D'_{\Til M}]$$
after push@-forward to $Z'\subset M$, as needed. This concludes the
proof of (3), in the equivalent formulation stated in
Claim~III.4.\qed\enddemo


\head \S4. Remarks and applications\endhead
\subhead \S4.1. $\mu$@-class and Parusi\'nski's Milnor number\endsubhead
Taking degrees in Theorem~I.5 gives
$$\int \cmp(X)=\int \frac{c(TM)}{c(\Cal L)}\cap[X]+ \int c(\Cal L)^{\dim
X}\cap (\mu_{\Cal L}(Y)^\vee\otimes_M \Cal L)$$ 
Now observe that the $\int$ picks up the term of degree $\dim M=\dim
X+1$ in the last term. Thinking of $c(\Cal L)^{\dim X}$ as $c(\Cal
L^{\oplus \dim X})$ and using (an immediate generalization of) \cite{A-F}:
$$\int c(\Cal L)^{\dim X}\cap (\mu_{\Cal L}(Y)^\vee\otimes \Cal L)
=\int c(\Cal L^{\oplus \dim X}\otimes \Cal L^\vee)\cap (\mu_{\Cal
L}(Y)^\vee\otimes \Cal L\otimes \Cal L^\vee)=\int \mu_{\Cal L}(Y)^\vee$$
Recalling that the degree of $\cmp(X)$ equals the Euler characteristic
of $X$, this proves:
\proclaim{Proposition IV.1} 
$$\int \mu_{\Cal L}(Y)=(-1)^{\dim M}\left(\chi(X)-\int
\frac{c(TM)}{c(\Cal L)}\cap[X]\right)\tag*$$
\endproclaim
Over $\Bbb C$, the right@-hand@-side in this formula equals
Parusi\'nski's generalization of Milnor's number (\cite{Parusi\'nski});
so this gives an alternative proof of Proposition~2.1 in
\cite{Aluffi1}, and extends to arbitrary fields of characteristic~0
the interpretation of $\mu_{\Cal L}(Y)$ as a measure of the difference
in the Euler characteristics of special vs.~general sections of a line
bundle (if $\Cal L$ has enough sections, the last term in (*) gives
$\chi(X_g)$ for a general section $X_g$ of $\Cal L$).

Conversely, at least if $\Cal L$ is ample enough, the formula in
Proposition~IV.1 suffices to prove Theorem~I `numerically', that is up to
taking degrees with respect to $\Cal L$. This is worked out in
\cite{Aluffi2}.

For isolated singularities on strong local complete intersections, a
statement analogous to Theorem~I.5 has been proved by T.~Suwa
\cite{Suwa}.

\subhead \S4.2. Blowing up $\mu$@-classes\endsubhead
The blow@-up formula proved in \S3 translates nicely in terms of
$\mu$@-classes.

Notations as in (3) from \S2: $Z\subset X\subset M$ is a nonsingular
subvariety of codimension $d$ in $M$ ($\dim M=n$), $\Til M=B\ell_Z M
@>>> \Til M$ denotes the blow@-up of $M$ along $Z$, and $X'$ denotes the
(scheme@-theoretic) inverse image of $X$ in $\Til M$. Also, $Y,Y'$ are
the singular schemes of $X,X'$ respectively. If $\Cal L$ denotes the
line bundle of $X$, note its pull@-back is the line@-bundle of $X'$.
Then the equality (3) we proved in \S 3:
$$\pi_*(c_*(X'))=c_*(X)+(d-1)\,c_*(Z)$$
(with $c_*$ as in \S1) becomes, in terms of $\mu$@-classes:
\proclaim{Proposition IV.2}
$$\pi_* \mu_{\Cal L}(Y')=\mu_{\Cal L}(Y)+(-1)^d(d-1) \mu_{\Cal L}(Z)$$
\endproclaim
\demo{Proof} Since $\Cal O(X')$ is the pull@-back of $\Cal O(X)$:
$$\align
\pi_*(c(T\Til M) &\cap s(X',\Til M))-c(TM)\cap s(X,M)\\
&=\pi_*\left(c(T\Til M)\cap\frac {\pi^*[X]}{1+\pi^* X}\right)
-c(TM)\cap\frac{[X]}{1+X}
\\ &=\left(\pi_*(c(T\Til M)\cap [\Til M])-c(TM)\cap [M]\right) \cdot
\frac{[X]}{1+X}
\endalign$$
Now we already observed in \S3.1 that
$$\pi_*(c(T\Til M)\cap [\Til M])-c(TM)\cap [M]=(d-1)\,c(TZ)\cap
[Z]\quad;$$
therefore (using the expression for $c_*$ in Theorem~I.5), (3) is
equivalent to:
$$\multline
\pi_*\left(c(\Cal L)^{n-1}\cap (\mu_{\Cal L}(Y')^\vee \otimes \Cal
L)\right) - c(\Cal L)^{n-1}\cap (\mu_{\Cal L}(Y)^\vee \otimes \Cal
L)\\
= (d-1)\, c(TZ)\cap [Z]-(d-1)\, c(TZ)\cap [Z]\cdot \frac{[X_1]}{1+X_1}
\endmultline$$
that is, to:
$$\pi_*\left(c(\Cal L)^{n-1}\cap (\mu_{\Cal L}(Y')^\vee \otimes \Cal
L)\right) - c(\Cal L)^{n-1}\cap (\mu_{\Cal L}(Y)^\vee \otimes \Cal
L)= (d-1)\, \frac {c(TZ)}{c(\Cal L)}\cap [Z]$$
Now we apply easy manipulations (see \cite{Aluffi2}, \S2):

---cap by $c(\Cal L)^{-(n-1)}$:
$$\pi_*\mu_{\Cal L}(Y')^\vee \otimes \Cal L - \mu_{\Cal L}(Y)^\vee
\otimes \Cal L = (d-1)\, \frac {c(TZ)}{c(\Cal L)^n}\cap [Z]$$

---$\otimes \Cal L^\vee$:
$$\pi_*\mu_{\Cal L}(Y')^\vee - \mu_{\Cal L}(Y)^\vee = (d-1)\,
c(\Cal L^\vee)^d \frac{c(TZ\otimes \Cal L^\vee)}{c(\Cal L\otimes \Cal
L^\vee)^n}\cap \frac {[Z]}{c(\Cal L^\vee)^d}$$

---clean up, dualize, apply Corollary~1.8 from \cite{Aluffi1}:
$$\pi_*\mu_{\Cal L}(Y') - \mu_{\Cal L}(Y) = (-1)^d(d-1)\,
c(T^*Z\otimes \Cal L)\cap [Z]=(-1)^d(d-1)\, \mu_{\Cal L}(Z)$$

This is the equality stated above, and we are done.\qed\enddemo

The relation of $\mu$@-classes stated above is of some independent
interest. It is related to a result in \cite{Parusi\`nski} (Lemma
2.2), which can be stated as the fact that the zero@-dimensional terms
of the sides have the same degree.

\subhead \S4.3. Contact of two hypersurfaces\endsubhead
Our proof of the main Theorem used very little of the good functoriality
properties of $\cmp(X)$: we proved just enough of them for $c_*(X)$ to
force this to equal $\cmp(X)$. After the fact, however, $c_*(X)$ inherits
the full set from $\cmp(X)$, and this reflects into facts about Segre
classes of singular schemes of hypersurfaces which we are not able to prove
otherwise, or which would require substantially more work by more
conventional techniques. While we plan to explore this elsewhere, we give a
few such examples in this subsection. 

Suppose $M_1$, $M_2$ are distinct nonsingular hypersurfaces of a
nonsingular ambient variety $M$. Then $X=M_1\cap M_2$ is a hypersurface of
both $M_1$ and $M_2$, with normal bundle $\Cal L_2=\Cal O(M_2)|_X$ in $M_1$
and $\Cal L_1=\Cal O(M_1)|_X$ in $M_2$. The singular scheme $Y$ of $X$ is
supported on the locus where $M_1$ and $M_2$ are {\it tangent:\/} we call
$Y$ the {\it contact scheme} of $M_1$ and $M_2$ in this case. What can be
said in general about $Y$?
\proclaim{Proposition IV.3} Under the above hypotheses (and with the
notation introduced in \S1.4),
$$s(Y,M_1)\otimes_{M_1} \Cal L_1= s(Y,M_2)\otimes_{M_2} \Cal L_2$$
\endproclaim
\demo{Proof} By Theorem I.4, we can compute $\cmp(X)$ by viewing it as a
hypersurface of $M_1$:
$$\cmp(X)=c(TM_1)\cap\left(s(X,M_1)+c(\Cal L_2)^{-1}\cap(s(Y,M_1)^\vee
\otimes_{M_1} \Cal L_2)\right)$$
or as a hypersurface of $M_2$:
$$\cmp(X)=c(TM_2)\cap\left(s(X,M_2)+c(\Cal L_1)^{-1}\cap(s(Y,M_2)^\vee
\otimes_{M_2} \Cal L_1)\right)$$
It follows that the right@-hand@-side of these expressions are equal. The
first summand in both is $c_F(X)$, so we get
$$\frac{c(TM_1)}{c(\Cal L_2)}\cap(s(Y,M_1)^\vee \otimes_{M_1} \Cal L_2)
=\frac{c(TM_2)}{c(\Cal L_1)}\cap(s(Y,M_2)^\vee \otimes_{M_2} \Cal L_1)$$
Capping with the inverse Chern class of the virtual tangent bundle of
$X$ and dualizing:
$$s(Y,M_1) \otimes_{M_1} \Cal L_2^\vee=s(Y,M_2) \otimes_{M_2} \Cal L_1^\vee$$
Tensoring both sides by $\Cal L_1\otimes\Cal L_2$ gives the statement of
the proposition.\qed\enddemo

We do not see any simple way to derive the result in Proposition~IV.3 more
directly. The result prompts us to define a class
$$\Cal S(Y,M)=s(Y,M_1)\otimes_{M_1} \Cal L_1\quad,$$
since we just showed that this is in a sense intrinsic to the contact
scheme and to the ambient variety. It would be interesting to study
properties of this class.

Next, observe that $M_1\cup M_2$ is a hypersurface of $M$, with line bundle
$\Cal O(M_1)\otimes\Cal O(M_2)$. The singular scheme $\overline X$ of
$M_1\cup M_2$ is supported on $X=M_1\cap M_2$, but `thicker' than $X$ along
$Y$. Now
$$\cmp(M_1\cup M_2)=\cmp(M_1)+\cmp(M_2)-\cmp(M_1\cap M_2)\tag**$$
(cf.~the beginning of \S~2). If the hypersurfaces involved are all divisors
with normal crossings, we proved this relation `by hand' in \S2 for the
class $c_*$ defined in \S1. As we have now proved that $c_*=\cmp$, we know
that this equality must hold regardless of how the hypersurfaces meet.
Using for example the expression for $c_*$ given in Theorem~I.4, this
gives a nontrivial relation among $s(\overline X,M)$ (on the
left@-hand@-side) and $s(Y,M_i)$ (on the right@-hand@-side).
Unraveling notations, the reader will check that this gives
\proclaim{Proposition IV.4} $s(\overline X,M)=s(X,M) + c(N_X M)^{-1}\cap \Cal
S(Y,M)\,$.\endproclaim

This is a sort of `residual intersection formula' (thinking of $Y$ as the
residual of $X$ in $\overline X$), and as such it could probably be proved
by judicious use of Proposition~9.2 in \cite{Fulton}, perhaps after blowing
up $M$ along $X$. The above argument seems however more direct at this
point, and the formula is perhaps simpler than it would be fair to
expect. Note: if $\overline X\subset M_1$, say, then standard residual
intersection formulas can be applied to $X\subset\overline X\subset M_1$
(as $X$ is a divisor of $M_1$), and do yield the formula stated in the
proposition. However, in general $\overline X$ is {\it not\/} contained in
either $M_1$ or $M_2$.

The above argument will work even if one of the hypersurfaces, say $M_2$,
is singular. In such a case $\overline X$ will be supported on $X=M_1\cap
M_2$ and on $W=$ the singular scheme of $M_2$. In terms of $c_*$, (**) then
says (all $\otimes$ in $M$ unless otherwise denoted)
$$\multline c(TM)\cap \left(\frac{[M_1]+[M_2]}{c(\Cal L_1\otimes\Cal L_2)}
+\frac 1{c(\Cal L_1\otimes\Cal L_2)}\cap (s(\overline X,M)^\vee \otimes\Cal
L_1 \otimes\Cal L_2)\right)\\
=c(TM)\cap \left(\frac{[M_1]}{c(\Cal L_1)}+\frac{[M_2]} {c(\Cal
L_2)}+\frac 1{c(\Cal L_2)}\cap (s(W,M)^\vee\otimes\Cal L_2)\right)\\
-c(T M_1)\cap \left(\frac{[M_1]\cdot [M_2]}{c(\Cal L_2)} -\frac
1{c(\Cal L_2)}(s(Y,M_1)^\vee \otimes_{M_1} \Cal L_2)\right)
\endmultline$$
(note: the last ${}^\vee$ is also taken in $M$, hence the change of sign in
the last $()$). The reader should have no difficulties (using
\cite{Aluffi2}) simplifying this expression to
\proclaim{Proposition IV.5} With the above notations,
$$s(\overline X,M)-s(X,M)=\frac 1{c(\Cal L_1)}\cap
\left(s(W,M)\otimes_M \Cal L_1\right) + \frac 1{c(\Cal L_2)}\cap
\left(s(Y,M_1)\otimes_M \Cal L_1\right)$$
\endproclaim

The residual intersection problem is in this case complicated enough that
we were not able to prove this relation otherwise.

\subhead \S4.4. A geometric application\endsubhead
If the singular scheme $Y$ of a hypersurface $X$ is nonsingular, then the  class
$c_*(X)$ of \S1 (and hence $\cmp(X)$) has a particularly simple form:
\proclaim{Proposition IV.6} Let $X$ be a hypersurface in a nonsingular variety
$M$, let $\Cal L=\Cal O(X)$ and assume that the singular scheme $Y$ of $X$
is nonsingular. Then
$$\cmp(X)=c_F(X)+(-1)^{\codim_M Y}\frac{c(TY)}{c(\Cal L)}\cap [Y]$$
\endproclaim
\demo{Proof} By \cite{Aluffi1}, Corollary~1.8, if $Y$ is nonsingular then
$$\mu_{\Cal L}(Y)=c(T^*Y\otimes\Cal L)\cap [Y]\quad;$$
hence by Theorem~I.5
$$\cmp(X)=c_F(X)+c(\Cal L)^{\dim X}\cap\left(\left(c(T^*Y\otimes\Cal L)\cap
[Y]\right)^\vee\otimes_M \Cal L\right)$$
Using \cite{Aluffi2}, \S2:
$$\align
c(\Cal L)^{\dim X}&\cap\left(\left(c(T^*Y\otimes\Cal L)\cap
[Y]\right)^\vee\otimes_M \Cal L\right) \\
&= (-1)^{\codim_M Y} c(\Cal L)^{\dim X}\left(\left(c(TY\otimes\Cal
L^\vee)\cap [Y]\right)\otimes_M \Cal L\right)\\
&=(-1)^{\codim_M Y} c(\Cal L)^{\dim X}\left(\frac{c(TY)}{c(\Cal L)^{\dim
Y}}\cap \frac{[Y]}{c(\Cal L)^{\codim_MY}}\right)\\
&=(-1)^{\codim_M Y} \frac{c(TY)}{c(\Cal L)}\cap [Y]
\endalign$$
which yields the statement.\qed\enddemo

As an application, consider again the situation at the beginning of \S4.3:
$M_1$, $M_2$ are nonsingular hypersurfaces with contact scheme $Y$
(=singular scheme of $M_1\cap M_2$). Proposition~IV.3 in \S4.3 spells out a
constraint imposed on the situation; in the particular case when $Y$ is
nonsingular, Proposition~IV.6 above allows us to rewrite this in a
particularly simple form:
\proclaim{Proposition IV.7} Assume $M_1$, $M_2$ are nonsingular hypersurfaces in
a nonsingular ambient variety, and let $Y$ be their contact scheme. Assume
$Y$ is nonsingular; then
$$M_1\cdot Y=M_2\cdot Y$$
\endproclaim
\demo{Proof} Let $X=M_1\cap M_2$. The last proposition can be used to
compute $\cmp(X)$ in two ways: considering $X$ as a hypersurface in $M_1$,
with normal bundle $\Cal L_2=\Cal O(M_2)|_X$, or as a hypersurface in
$M_2$, with normal bundle $\Cal L_1=\Cal O(M_1)|_X$. This gives:
$$c_F(X)+(-1)^{\codim_M Y}\frac{c(TY)}{c(\Cal L_2)}\cap [Y]=
c_F(X)+(-1)^{\codim_M Y}\frac{c(TY)}{c(\Cal L_1)}\cap [Y]$$
and hence
$$c(\Cal L_1)\cap [Y]=c(\Cal L_2)\cap [Y]$$
from which the stated formula follows.\qed\enddemo

As a concrete application, say the ambient variety is a projective space,
$d_i=\deg M_i$, and the contact scheme $Y$ of $M_1$, $M_2$ is nonsingular
and positive dimensional; then the statement is that necessarily
$d_1=d_2$. It is easy to produce examples of hypersurfaces of the same
degree and having nonsingular contact {\it scheme:} for instance, the quadrics
$$\align
x^2+y^2+z^2+w^2&=0\\
x^2+y^2+z^2+2 w^2&=0
\endalign$$
in $\P^3$ meet along a double conic, so have contact scheme equal to
a nonsingular plane conic. Proposition~IV.7 shows that no such example
can be concocted with smooth hypersurfaces of {\it different} degrees.



\enddocument